\documentclass[ps,preprint,dvips,twoside]{imsart}

\RequirePackage[colorlinks,citecolor=blue,urlcolor=blue]{hyperref}
\usepackage{url,breakurl}
\usepackage{amssymb,amsmath}
\usepackage{graphicx}

\newtheorem{theorem}{Theorem}[section]
\newtheorem{lemma}[theorem]{Lemma}
\newtheorem{proposition}[theorem]{Proposition}
\newtheorem{definition}[theorem]{Definition}
\newtheorem{corollary}[theorem]{Corollary}
\newtheorem{conjecture}[theorem]{Conjecture}

\doi{10.1214/13-PS214}
\pubyear{2014}
\volume{11}
\issue{0}
\firstpage{121}
\lastpage{160}

\allowdisplaybreaks

\begin{document}

\begin{frontmatter}
\title{Statistical properties of zeta functions'~zeros}
\runtitle{Statistics of zetas' zeros}

\date{May 2013}
\begin{aug}
\author{\fnms{Vladislav} \snm{Kargin}\corref{}\ead
[label=e1]{vladislav.kargin@gmail.com}}
\address{282 Mosher Way, Palo Alto, CA 94304\\
\printead{e1}}

\runauthor{V. Kargin}
\end{aug}

\begin{abstract}
The paper reviews existing results about the statistical distribution of
zeros for three main types of zeta functions: number-theoretical,
geometrical, and dynamical. It provides necessary background and
some details about the proofs of the main results.
\end{abstract}



% history:
\received{\smonth{2} \syear{2013}}
\end{frontmatter}

\tableofcontents

%s1 ###
\section{Introduction}

The distribution of zeros of Riemann's zeta function is one of the central
problems in modern mathematics. The famous Riemann conjecture states that
all of these zeros are on the critical line $\mathrm{Re}s=1/2$, and the
accumulated numerical evidence supports this conjecture as well as a more
precise statement that these zeros behave like eigenvalues of large random
Hermitian matrices. While these statements are still conjectural, a great
deal is known about the statistical properties of Riemann's zeros and the
zeros of closely related functions. In this report we aim to summarize
findings in this research area.

We give necessary background information, and we cover the three main types
of zeta functions: number-theoretical, Selberg-type, and dynamical zeta
functions.

Some interesting and important topics are left outside of the scope of this
report. For example, we do not discuss quantum arithmetic chaos or
characteristic polynomials of random matrices.

The paper is divided in three main sections according to the type of the
zeta function we discuss. Inside each section we tried to separate the
discussion of the properties of zeros at the global and local scales.

Let us briefly describe these types of zeta functions and their
relationships. First, the number-theoretical zeta functions come from
integers in number fields and their generalizations. Due to the
additive and
multiplicative structures of the integers, and in particular due to the
unique decomposition in prime factors, the zeta functions have the Euler
product formula $\zeta(s) =\prod_{\mathfrak{p}}(1-(N\mathfrak{p})^{s})^{-1}$
and a functional equation, $\zeta( 1-s) =c(s) \zeta(s)$, with the
multiplier $c(s)$ equal to a ratio of Gamma functions.

It is an important discovery of Hecke that one can define number-theoretical
zeta functions in a different, and potentially more general way if one
starts with modular forms, which are functions on the space of
two-dimensional lattices that are invariant relative to the change of scale.
If they are considered as functions of the basis $( z,1) $ then
they become functions of $z$ invariant relative to an action of the
group $%
SL_{2}( \mathbb{Z}) $. They are periodic and hence can be written
as $\sum_{n\geq0}c_{n}\exp( 2\pi inz),$ where $z$ is the ratio
of the periods of the lattice. If $c_{0}=0,$ then one can define a zeta
function $\sum c_{n}n^{-s}$. It turns out that this zeta function satisfies
a functional equation. In addition, if the original modular form is an
eigenvector for certain operators (the Hecke operators), then the zeta
function will have the Euler product property.

It appears that all number-theoretical zeta functions can be obtained by
using this construction. However, the class of modular zeta functions is
wider. Indeed, zeta functions with an arithmetic flavor can come not only from number
fields but also from other algebraic objects, for example, from elliptic
curves. A significant development occurred recently when it was proved that
all zeta functions associated with elliptic curves come from modular zeta
functions. Among other applications, this discovery was a key to the proof
of Fermat's last theorem. (The existence of a non-trivial solution for $
x^{n}+y^{n}=z^{n},$ $n>2,$ would imply the existence of a non-modular zeta
function for an elliptic curve.)

An important representative of the second class of zetas is Selberg's zeta
function, which is essentially a generating function for the lengths of
closed geodesics on a surface with constant negative curvature. In more
detail, let $\mathbb{H}$ be the upper half-plane with the hyperbolic metric
and $\Gamma$ be a discrete subgroup of $SL_{2}( \mathbb{Z})$.
Selberg showed that certain sums over eigenvalues of the Laplace
operator on
Riemann surface $\Gamma\backslash\mathbb{H}$ can be related to sums over
closed geodesics of $\Gamma\backslash\mathbb{H}.$ This relation is called
Selberg's trace formula. Selberg's zeta function is constructed in such a
way that it has the same relation to Selberg's trace formula as Riemann's
zeta function has to the so-called ``explicit
formula'' for sums over Riemann's zeros.

While Selberg's zeta function resembles Riemann's zeta in some features,
there are significant differences. In particular, the statistical behavior
of its zeros depends on the group $\Gamma$ and often it is significantly
different from the behavior of Riemann's zeros.

The third class of zetas, the dynamical zeta functions, are generating
functions for the lengths of closed orbits of a map $f$ that sends a
set $M$
to itself. The most spectacular example of these functions is Weil's zeta
functions of algebraic varieties over finite fields, which can be described as follows.

Let $M$ be the algebraic closure of an algebraic variety embedded in $%
\mathbb{F}_{q}^{n}$ (where $\mathbb{F}_{q}$ is a finite field), and let $f$
be the Frobenius map: $( x_{1},\ldots,x_{n}) \rightarrow$ $%
( x_{1}^{q},\ldots,x_{n}^{q})$. Then, the dynamic zeta
associated to $f$ is called Weil's zeta function. These functions are
remarkable since the Riemann conjecture is proved for them: It is known that
their zeros are located on a circle that corresponds to the line
$\mathrm{Re}%
z=1/2.$ Moreover, much is known about the statistical distribution of these
zeros.\looseness=-1

A particular case of these zeta functions, Weil's zeta functions for curves,
can be understood as number-theoretic zeta functions for finite extensions
of the field $\mathbb{F}_{q}( x)$. Hence, Weil's  zeta functions work as a bridge
between dynamic and number-theoretic zeta functions. In addition, Selberg's
zeta function can be understood as the dynamic zeta function for the
geodesic flow on the surface $\Gamma\backslash\mathbb{H}$. This shows that
the three classes of zeta functions are intimately related to each
other.\looseness=-1

With this overview in mind, we now come to a more detailed description of
available results about the statistics of zeta function zeros.

%s2 ###
\section{Number-theoretical zetas}

%s2.1 ###
\subsection{Riemann's zeta}

There are several excellent sources on Riemann's zeta and Dirichlet
L-functions, for example the books by Davenport \cite{davenport67} and
Titchmarsh \cite{titchmarsh86}. In addition, a very good reference for all
topics in this report is provided by Iwaniec and Kowalski's book \cite%
{iwaniec_kowalski04}.

By definition Riemann's zeta function is given by the series
%
%e1 ###
\begin{equation}
\zeta\left( s\right) =\sum_{n=1}^{\infty}\frac{1}{n^{s}},
\label{definition_Riemann_zeta}
\end{equation}
for $\mathrm{Re}s>1.$ It can be analytically continued to a meromorphic
function in the entire complex plane and it satisfies the functional
equation
%
%e2 ###
\begin{equation}
\pi^{-\frac{s}{2}}\Gamma\left( \frac{1}{2}s\right) \zeta\left(
s\right)
=\pi^{-\frac{1-s}{2}}\Gamma\left( \frac{1}{2}\left( 1-s\right) \right)
\zeta\left( 1-s\right). \label{functional_equation_Riemann}
\end{equation}
Indeed, we can relate $\zeta( s) $ to the series $\theta
(x) =\sum_{n=1}^{\infty}e^{-n^{2}\pi x}:$
\begin{equation*}
\frac{\Gamma\left( \frac{1}{2}s\right) \zeta\left( s\right) }{\pi
^{s/2}}%
=\sum_{n=1}^{\infty}\int_{0}^{\infty}x^{s/2-1}e^{-n^{2}\pi
x}dx=\int_{0}^{\infty}x^{s/2-1}\theta\left( x\right) dx.
\end{equation*}
(As an aside remark, this representation can be used as a starting
point for
some surprising connections of the Riemann zeta function with the Brownian
motion and Bessel processes, see \cite{biane_pitman_yor01}.) From the
identities for the Jacobi theta-functions, implied by the Poisson summation
formula, it follows that
%
%e3 ###
\begin{equation}
2\theta\left( x\right) +1=\frac{1}{\sqrt{x}}\left( 2\theta\left( \frac
{1}{%
\sqrt{x}}\right) +1\right). \label{modularity_condition}
\end{equation}
Writing
\begin{equation*}
\int_{0}^{\infty}x^{s/2-1}\theta\left( x\right)
dx=\int_{0}^{1}x^{s/2-1}\theta\left( x\right) dx+\int_{1}^{\infty
}x^{s/2-1}\theta\left( x\right) dx,
\end{equation*}
and applying the identity to the integral from $0$ to $1,$ we find that
%
%e4 ###
\begin{equation}
\frac{\Gamma\left( \frac{1}{2}s\right) \zeta\left( s\right) }{\pi
^{s/2}}=%
\frac{1}{s\left( s-1\right) }+\int_{1}^{\infty}\left(
x^{-s/2-1/2}+x^{s/2-1}\right) \theta\left( x\right) dx,
\label{theta_representation}
\end{equation}
which is symmetric relative to the change $s\rightarrow1-s.$

In addition, formula (\ref{theta_representation}) implies that the
function $%
s( s-1) \pi^{-s/2}\Gamma( \frac{1}{2}s) \zeta(
s) $ is entire. Hence, the only pole of $\zeta(s)$ is at $s=1,$ and $%
\zeta( s) $ has zeros at $-2,$ $-4,$ \ldots, that correspond to
poles of $\Gamma( \frac{1}{2}s).$ These zeros are called \emph{%
trivial}, while all others are called \emph{non-trivial}. We will order the
non-trivial zeros according to their imaginary part and denote them by
$\rho
_{k}.$ By the functional equation, $\rho_{k}$ are located symmetrically
relative to the line $\mathrm{Re}s=1/2,$ which is called \emph{critical},
and it is known that $0<\mathrm{Re}\rho_{k}<1.$ The Riemann's hypothesis
asserts that all non-trivial zeros are on the critical line.

The second fundamental property of Riemann's zeta is the Euler product
formula:
%
%e5 ###
\begin{equation}
\zeta\left( s\right) =\prod_{p}\left( 1-\frac{1}{p^{s}}\right) ^{-1}
\label{product_formula_Riemann}
\end{equation}
valid for $\mathrm{Re}s>1.$ It follows from the existence and
uniqueness of
prime factorization. This formula is a starting point for a very important
idea which relates the sums over prime numbers and sums over zeta function
zeros. For example, the Riemann-von Mangoldt formula says that
%
%e6 ###
\begin{equation}
\sum_{n\leq x}\Lambda\left( n\right) =x-\sum_{\rho}\frac{x^{\rho
}}{\rho}%
+\sum_{n}\frac{x^{-2n}}{2n}-\frac{\zeta^{\prime}}{\zeta}\left(
0\right),
\label{Riemann_Mangoldt}
\end{equation}
where $\Lambda( n) =\log p,$ if $n$ is a prime $p$ or a power of
$p,$ and otherwise $\Lambda( n) =0$. (The order of summation
over zeros can be important here, so it is assumed that in computing $%
\sum_{\rho}\frac{x^{\rho}}{\rho},$ one takes all zeros with imaginary
part between $-T$ and $T$, and then let $T\rightarrow\infty.$ In addition,
if $x$ is a prime power the formula has to be modified by subtracting
$\frac{%
1}{2}\Lambda( x) $ on the left-hand side.)

The idea of the proof of (\ref{Riemann_Mangoldt}) is to take the logarithmic
derivative of (\ref{product_formula_Riemann}):
%
%e7 ###
\begin{equation}
\frac{\zeta^{\prime}}{\zeta}\left( s\right) =-\sum_{n=2}^{\infty
}\frac{%
\Lambda\left( n\right) }{n^{s}}, \label{log_derivative_Riemann}
\end{equation}
and then to use the following formula with $y=x/n$ and $c>0$:
\begin{equation*}
\frac{1}{2\pi i}\int_{c-i\infty}^{c+i\infty}y^{s}\frac{ds}{s}=\left\{
\begin{array}{cc}
0 & \text{if }0<y<1, \\
\frac{1}{2} & \text{if }y=1, \\
1 & \text{if }y>1.%
\end{array}
\right.
\end{equation*}
 One integrates (\ref%
{log_derivative_Riemann}) against the test function $x^s/s,$ in order to pick out the terms in the series in (\ref%
{log_derivative_Riemann}) with $n\leq x.$

From (\ref{log_derivative_Riemann}), one gets
\begin{equation*}
\sum_{n\leq x}\Lambda\left( n\right) =\frac{1}{2\pi i}\int_{c-i\infty
}^{c+i\infty}\left[ -\frac{\zeta^{\prime}}{\zeta}\left( s\right)
\right]
x^{s}\frac{ds}{s}.
\end{equation*}
Moving the line of integration away to infinity on the left and collecting
the residues at the poles one finds formula (\ref{Riemann_Mangoldt}). (See
Chapter 17 in Davenport \cite{davenport67} for a detailed proof.)

By a similar method one can obtain:
%
%e8 ###
\begin{equation}
\sum_{n\leq x}\frac{\Lambda\left( n\right) }{n^{s}}=\frac{x^{1-s}}{1-s}
-\sum_{\rho}\frac{x^{\rho-s}}{\rho-s}+\sum_{n}\frac{x^{-2n-s}}{2n+s}-
\frac{\zeta^{\prime}}{\zeta}\left( s\right). \label{Landau}
\end{equation}
While formula (\ref{Riemann_Mangoldt}) is useful to study the distribution
of primes if something is known about the distribution of zeros,
formula (%
\ref{Landau}) can be used in the reverse direction to study the
behavior of $%
\frac{\zeta^{\prime}}{\zeta}( s) $ if some information is
known about primes.

Selberg discovered a variant of this formula that avoids the problem of
conditional convergence in the sum over the zeta zeros. Define
\begin{equation*}
\Lambda_{x}\left( n\right) =\left\{
\begin{array}{cc}
\Lambda\left( n\right), & \text{ for }1\leq n\leq x, \\
\Lambda\left( n\right) \left( \frac{\log^{2}\frac{x^{3}}{n}-2\log^{2}
\frac{x^{2}}{n}}{2\log^{2}x}\right), & \text{for }x\leq n\leq x^{2},
\\
\Lambda\left( n\right) \frac{\log^{2}\frac{x^{3}}{n}}{2\log^{2}x}, &
\text{for }x^{2}\leq n\leq x^{3}.%
\end{array}
\right.
\end{equation*}
Then, (Lemma 10 in \cite{selberg46a} )
%
%e9 ###
\begin{eqnarray}
\frac{\zeta^{\prime}}{\zeta}\left( s\right) &=&-\sum_{n\leq
x^{3}}\frac{%
\Lambda_{x}\left( n\right) }{n^{s}}+\frac{1}{\log^{2}x}\sum_{\rho
}\frac{%
x^{\rho-s}\left( 1-x^{\rho-s}\right) ^{2}}{\left( s-\rho\right) ^{3}}
\label{explicit_formula_Selberg} \\
&&{}+\frac{x^{1-s}\left( 1-x^{1-s}\right) ^{2}}{\log^{2}x\cdot\left(
1-s\right) ^{3}}+\frac{1}{\log^{2}x}\sum_{q=1}^{\infty}\frac{%
x^{-2q-s}\left( 1-x^{-2q-s}\right) ^{2}}{\left( 2q+s\right) ^{3}}.
\notag
\end{eqnarray}

Another useful variant is sometimes called Weil's and sometimes Delsarte's
explicit formula (see \cite{conrey03} and \cite{moreno05}). Suppose
that $%
H( s) $ is an analytic function in the strip $-c\leq\mathrm{Im}$
$s\leq1+c$ (for $c>0$) and that $\vert H( \sigma+it)
\vert\leq A( 1+\vert t\vert) ^{-(
1+\delta) }$ uniformly in $\sigma$ in the strip. Let $h(
t) =H( \frac{1}{2}+it) $ and define
\begin{equation*}
\widehat{h}\left( x\right) =\int_{\mathbb{R}}h\left( t\right) e^{-2\pi
itx}dt.
\end{equation*}
(Note that analyticity of $H( s) $ implies that $\widehat{h}%
( x) $ has finite support.) Then,
%
%e10 ###
\begin{eqnarray}
\sum_{\rho}H\left( \rho\right) &=&H\left( 0\right) +H\left( 1\right) -
\frac{1}{2\pi}\sum_{n=1}^{\infty}\frac{\Lambda\left( n\right) }{\sqrt
{n}}%
\left[ \widehat{h}\left( \frac{\log n}{2\pi}\right) +\widehat{h}\left(
-%
\frac{\log n}{2\pi}\right) \right] \notag\\
&&{}-\frac{1}{2\pi}\int_{-\infty}^{\infty}h\left( t\right) \Psi\left(
t\right) dt, \label{formula_delsarte}
\end{eqnarray}
where
\begin{equation*}
\Psi\left( t\right) =\frac{\Gamma^{\prime}}{\Gamma}\left( \frac
{1}{2}%
+it\right) +\frac{\Gamma^{\prime}}{\Gamma}\left( \frac{1}{2}-it\right
).
\end{equation*}

The idea of the proof is similar to the proof of the Riemann-von Mangoldt
formula. One starts with the formula
\begin{eqnarray*}
\frac{1}{2\pi i}\int_{2-i\infty}^{2+i\infty}\left[ -\frac{\zeta
^{\prime}%
}{\zeta}\left( s\right) \right] H\left( s\right) ds &=&-\sum
_{n=1}^{\infty}%
\frac{\Lambda\left( n\right) }{n^{1/2}}\frac{1}{2\pi}\int_{-\infty
}^{\infty}h\left( t\right) e^{-i\left( \log n\right) t}dt \\
&=&-\frac{1}{2\pi}\sum_{n=1}^{\infty}\frac{\Lambda\left( n\right) }{%
n^{1/2}}\widehat{h}\left( \frac{\log n}{2\pi}\right).
\end{eqnarray*}

Next, one can move the line of integration to $( -1-i\infty
,-1+i\infty) $ and use the calculus of residues and the functional equation to
obtain formula (\ref{formula_delsarte}).

%f1 ###
\begin{figure}[t]
\includegraphics[width=11cm]{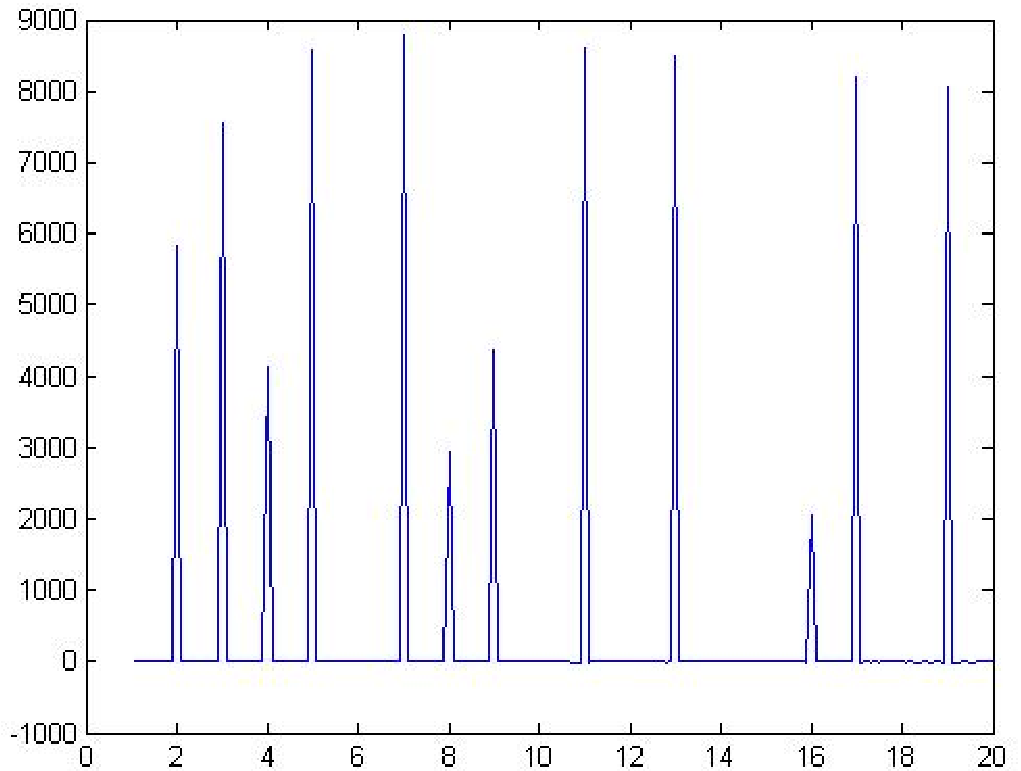}
\caption{Sums of $\cos(\protect\gamma\log(x))$ over the first $100{,}000$
zeros ($T\approx75\times10^{3}$). The horizontal axis shows $x$. }\label{fig:Landau}
\end{figure}

A wonderful illustration for the significance of explicit formulas is given
by Landau's formula:
\begin{equation*}
\sum_{0<\mathrm{Im}\rho<T}x^{\rho}=-\frac{T}{2\pi}\Lambda\left(
\sqrt{x}%
\right) +O\left( \log T\right).
\end{equation*}
Assuming Riemann's hypothesis, we write $\rho=\frac{1}{2}+i\gamma,$ and then
\begin{equation*}
\sum_{0<\gamma<T}\cos\left( \gamma\log x\right) =-\frac{T}{2\pi
}\frac{%
\Lambda\left( x\right) }{x}+O\left( \log T\right).
\end{equation*}
Note that the right hand side has a spike when $x$ is a prime power.
This is
illustrated in Figure \ref{fig:Landau}.\vadjust{\goodbreak}

The most general explicit formula was derived by Weil (see \cite{weil52}).
We will not present it here since it involves adelic language and this would
take us too far afield.\vspace*{-2pt}

%s2.1.1 ###
\subsubsection{Statistics of zeros on global scale}\vspace*{-2pt}

Let $\mathcal{N}(T)$ denote the number of zeros with the imaginary part
strictly between $0$ and $T$. If there is a zero with imaginary part equal
to $T$, then we count this zero as $1/2$. Define
\begin{equation*}
S(T):=\frac{1}{\pi}\mathrm{Im}\log\zeta\left(\frac{1}{2}+iT\right),
\end{equation*}
where the logarithm is calculated by continuous variation along the
contour $%
\sigma+iT,$ with $\sigma$ changing from $+\infty$ to $1/2.$

By applying the argument principle to $\zeta$ and utilizing the functional
equation, it is possible to show (see Chapter 15 in \cite{davenport67}) that
\begin{equation*}
\mathcal{N}(T)=\frac{T}{2\pi}\log\frac{T}{2\pi e}+\frac{7}{8}+S\left(
T\right) +O\left( \frac{1}{1+T}\right).
\end{equation*}
Let
\begin{equation*}
X\left( t\right) :=\frac{\sqrt{2}\pi S(t)}{\sqrt{\log\log t}}.
\end{equation*}
Then, we have the following theorem by Selberg. (See Theorem 6 in \cite%
{selberg46a}.)\vadjust{\goodbreak}

\begin{theorem}[Selberg]
\label{theorem_selberg} Let $T^{a}\leq H\leq T^{2},$ where $\frac{1}{2}%
<a\leq1.$ Then for every $k\geq1$ \
\begin{equation*}
\frac{1}{H}\int_{T}^{T+H}\left\vert X\left( t\right) \right\vert
^{2k}dt=%
\frac{2k!}{k!2^{k}}+O(\left( \log\log T\right) ^{-1/2}),
\end{equation*}
with the constant in the remainder term depending on $k$ and $a$ only.
\end{theorem}

In other words, if $t$ is chosen randomly in the interval $[0,T]$ and $%
T\rightarrow\infty,$ then $X( t) $ converges in distribution to
a Gaussian random variable. Note that the Riemann Hypothesis is not assumed
in this result. Under the assumption of the Riemann Hypothesis, Selberg
proved an analogous result with a better error term assuming only $a>0$
(Theorem 3 in \cite{selberg44}).

Recently, Selberg's result was generalized by Bourgade, who determined the
correlation structure of $X( t) $ on relatively small scales. The
basis for this development is the following result (cf. Theorem 1.1 in
\cite%
{bourgade10}). Take a large $t>0$ and small $\varepsilon_{t}\geq0$ and
look at $l$ points on the line $\mathrm{Re}z=\frac{1}{2}+\varepsilon_{t}$
with imaginary parts $\omega t+f_{t}^{( i) },$ $1\leq i\leq l,$
where $\omega$ is a random variable uniformly distributed on $[ 0,1%
].$ Note that the randomness $\omega$ is the same for all points so
the distance between the points is measured by offsets $f_{t}^{(
i) }.$ We assume that $\varepsilon_{t}\rightarrow0,$ so the points
approach the critical line as $t\rightarrow\infty$. The first case is
when $%
\varepsilon_{t}\gg1/\log t$ and we look at the scale $\vert
f_{t}^{( j) }-f_{t}^{( i) }\vert\approx
\varepsilon_{t}^{c_{ij}},$ with $c_{ij}\geq0.$ The second case is
when $%
\varepsilon_{t}\ll1/\log t,$ for example, $\varepsilon_{t}=0.$ In this
case, we look at the scale $\vert f_{t}^{( j)
}-f_{t}^{( i) }\vert\approx( 1/\log t)
^{c_{ij}}.$ It turns out that after a proper normalization the values
of the logarithm of the Riemann zeta function at this cluster of points converge to a non-trivial
multivariate Gaussian distribution.

\begin{theorem}[Bourgade]
\label{theorem_bourgade} Let $\omega$ be uniform on $(0,1)$, $\epsilon
_{t}\rightarrow0$, $t\rightarrow\infty,$ $\epsilon_{t}\gg1/\log
t$, and
functions $0\leq f_{t}^{( 1) }<\cdots<f_{t}^{( l)
}<c<\infty.$ Suppose that for all $i\neq j,$
%
%e11 ###
\begin{equation}
\frac{\log\left\vert f_{t}^{\left( j\right) }-f_{t}^{\left( i\right)
}\right\vert}{\log\epsilon_{t}}\rightarrow c_{i,j}\in\left[
0,\infty%
\right]. \label{bourgade_scale1}
\end{equation}
Then the vector
%
%e12 ###
\begin{equation}
\frac{1}{\sqrt{-\log\epsilon_{t}}}\left( \log\zeta\left( \frac{1}{2}
+\varepsilon_{t}+if_{t}^{\left( 1\right) }+i\omega t,\ldots,\frac
{1}{2}%
+\varepsilon_{t}+if_{t}^{\left( l\right) }+i\omega t\right) \right)
\label{bourgade_scale2}
\end{equation}
converges in law to a complex Gaussian vector $( Y_{1},\ldots
,Y_{l}) $ with the zero mean and covariance function
\begin{equation*}
\mathrm{Cov}\left( Y_{i},Y_{j}\right) =\left\{
\begin{array}{cc}
1 & \text{if }i=j, \\
1\wedge c_{i,j} & \text{if }i\neq j.%
\end{array}
\right.
\end{equation*}
Moreover, the above result remains true if $\epsilon_{t}\ll1/\log t$,
replacing the normalization $-\log\varepsilon_{t}$ with $\log\log t$
in (%
\ref{bourgade_scale1}) and (\ref{bourgade_scale2}).
\end{theorem}

This theorem implies the following result for the zeros of the Riemann zeta
function (cf. Corollary 1.3 in \cite{bourgade10}). Let
\begin{equation*}
\Delta\left( t_{1},t_{2}\right) =\left( \mathcal{N}(t_{2})-\frac
{t_{2}}{%
2\pi}\log\frac{t_{2}}{2\pi e}\right) -\left( \mathcal{N}(t_{1})-\frac
{t_{1}%
}{2\pi}\log\frac{t_{1}}{2\pi e}\right),
\end{equation*}
which represents the number of zeros with the imaginary part between
$t_{1}$ and
$t_{2}$ minus its deterministic prediction. Then the claim is that this
excess number of zeros in an interval with the length of order $( \log
t) ^{-\delta}$ ($0<\delta<1)$ is a Gaussian variable with the
variance proportional to $( 1-\delta) \log\log t.$ Moreover,
the limiting Gaussian process has an interesting covariance structure.

\begin{corollary}[Bourgade]
\label{corollary_bourgade} Let $K_{t}$ be such that, for some
$\varepsilon
>0 $ and all $t$, $K_{t}>\varepsilon$. Suppose $\log K_{t}/\log\log
t\rightarrow\delta\in\lbrack0,1),$ as $t\rightarrow\infty.$ Then the
finite-dimensional distributions of the process
\begin{equation*}
\frac{\Delta\left( \omega t+\alpha/K_{t},\omega t+\beta/K_{t}\right)
}{%
\frac{1}{\pi}\sqrt{\left( 1-\delta\right) \log\log t}},\text{ }0\leq
\alpha<\beta<\infty,
\end{equation*}
converge to those of a centered Gaussian process ($\widetilde{\Delta}(
\alpha,\beta) $, $0\leq\alpha<\beta<\infty$) with the covariance
structure
\begin{equation*}
\mathbb{E}\left( \widetilde{\Delta}\left( \alpha,\beta\right)
\widetilde{%
\Delta}\left( \alpha^{\prime},\beta^{\prime}\right) \right) =\left
\{
\begin{array}{cc}
1 & \text{if }\alpha=\alpha^{\prime},\text{ and }\beta=\beta
^{\prime},
\\
1/2 & \text{if }\alpha=\alpha^{\prime},\text{ and }\beta\neq\beta
^{\prime}, \\
1/2 & \text{if }\alpha\neq\alpha^{\prime},\text{ and }\beta=\beta
^{\prime}, \\
-1/2 & \text{if }\beta=\alpha^{\prime}, \\
0 & \text{elsewhere.}%
\end{array}
\right.
\end{equation*}
\end{corollary}

Note that since the average spacing between zeros is $1/\log t,$ hence the
number of zeros in the interval $( \omega t+\alpha/K_{t},\omega
t+\beta/K_{t}) $ is of order $( \log t) ^{1-\delta}.$

This result perfectly corresponds to a result of Diaconis and Evans about
eigenvalue fluctuations of random unitary matrices (cf. Theorem 6.3 in
\cite%
{diaconis_evans01}).

The key to both Selberg and Bourgade's results is Selberg's approximation
for the function $S( t) $ (cf. Theorem 4 in \cite{selberg46a}).

\begin{proposition}[Selberg]
\label{proposition_approximation}Suppose $k\in\mathbb{Z}^{+},$ $0<a<1.$
Then there exists $c_{a,k}>0$ such that for any $1/2\leq\sigma\leq1$
and $%
t^{a/k}\leq x\leq t^{1/k},$ it is true that
\begin{equation*}
\frac{1}{t}\int_{1}^{t}\left\vert\log\zeta\left( \sigma+is\right)
-\sum_{p\leq x^{3}}\frac{1}{p^{\sigma+is}}\right\vert^{2k}ds\leq c_{a,k}.
\end{equation*}
\end{proposition}

The proof of this statement, in turn, depends on Selberg's formula (\ref
{explicit_formula_Selberg}).

%s2.1.2 ###
\subsubsection{Statistics of zeros on local scale}

Now, assume the Riemann Hypothesis and suppose that we are interested in
calculating local statistics for the pairs of the zeta zeroes, for example,
in $\sum_{0<\gamma,\gamma^{\prime}\leq T}r( ( \gamma-\gamma
^{\prime}) \frac{\log T}{2\pi}),$\vspace*{1pt} where $r( x) $
is a test function. By representing $r( x) $ as the Fourier
transform, we can write this statistic differently,
%
%e13 ###
\begin{equation}
\sum_{0<\gamma,\gamma^{\prime}\leq T}r\left( \left( \gamma-\gamma
^{\prime}\right) \frac{\log T}{2\pi}\right) w\left( \gamma-\gamma
^{\prime}\right) =\frac{T}{2\pi}\log T\int_{-\infty}^{\infty}F\left(
\alpha\right) \widehat{r}\left( \alpha\right) d\alpha,
\label{pair_correlations}
\end{equation}
where we added a convenient weighting function $w(u)=4/( 4+u^{2})
.$ Here $\widehat{r}( \alpha) =\int_{-\infty}^{\infty}r(
u) e^{-2\pi i\alpha u}du,$ and
\begin{equation*}
F\left( \alpha\right) =F\left( \alpha,T\right) :=\frac{2\pi}{T\log T}
\sum_{0<\gamma,\gamma^{\prime}<T}T^{i\alpha\left( \gamma-\gamma
^{\prime}\right) }w\left( \gamma-\gamma^{\prime}\right).
\end{equation*}
Consequently, all information about local statistics is encoded in the
function $F( \alpha).$ Montgomery proved the following result
\cite{montgomery73}. (In fact, the estimate holds uniformly throughout $
0\leq\alpha\leq1,$ as was later proved by Goldston in \cite{goldston81}.)

\begin{theorem}[Montgomery]
\label{theorem_Montgomery} For real $\alpha,$ the function $F( \alpha
) $ is real and\break  $F( \alpha) =F( -\alpha).$ If
$T>T_{0}( \varepsilon),$ then $F( \alpha) \geq
-\varepsilon$ for all $\alpha.$ For fixed $\alpha$ satisfying $0\leq
\alpha<1,$ we have
\begin{equation*}
F\left( \alpha\right) =\left( 1+o\left( 1\right) \right) T^{-2\alpha
}\log
T+\alpha+o\left( 1\right),
\end{equation*}
as $T$ tends to infinity; this holds uniformly for $0\leq\alpha\leq
1-\varepsilon.$
\end{theorem}

This result allows us to calculate the local statistics for smooth test
functions that have their Fourier transforms compactly supported on the
interval $[ -1,1].$

Montgomery conjectures that
\begin{equation*}
F\left( \alpha\right) =1+o\left( 1\right)
\end{equation*}
for $\alpha\geq1,$ uniformly in bounded intervals. If it is true,
then it
can be used to calculate the local statistics for a much larger class of
test functions. In particular, Montgomery's conjecture can also be
formulated as follows.

\begin{conjecture}[Montgomery]
\label{conjecture_Montgomery}For fixed $\alpha<\beta,$
\begin{equation*}
\sum_{\substack{ 0<\gamma,\gamma^{\prime}<T \\ 2\pi\alpha/\log
T<\left| \gamma-\gamma^{\prime}\right| <2\pi\beta/\log T}}1\sim
\left(
\int_{\alpha}^{\beta}\left[ 1-\left( \frac{\sin\pi u}{\pi u}\right)
^{2}%
\right] du+\delta\left( \alpha,\beta\right) \right) \frac{T}{2\pi
}\log T,
\end{equation*}
as $T$ goes to infinity. Here $\delta( \alpha,\beta) =1$ if $%
0\in[ \alpha,\beta],$ $\delta( \alpha,\beta) =0$
otherwise.
\end{conjecture}

The proof of Montgomery's theorem is based on the analysis of the following
variant of the explicit formula.

\begin{lemma}[Montgomery]
\label{lemma_montgomery}If $1<\sigma<2$ and $x\geq1$ then
\begin{align*}
\sum_{\gamma}\frac{(2\sigma-1)x^{i\gamma}}{\left( \sigma-%
\frac{1}{2}\right) +\left( t-\gamma\right) ^{2}} &=-x^{-1/2}\left(
\sum_{n\leq x}\Lambda\left( n\right) \left( \frac{x}{n}\right)
^{1-\sigma
+it}+\sum_{n>x}\Lambda\left( n\right) \left( \frac{x}{n}\right)
^{\sigma
+it}\right) \\
&\quad +x^{1/2-\sigma+it}\left( \log\tau+O_{\sigma}\left( 1\right) \right)
+O_{\sigma}\left( x^{1/2}\tau^{-1}\right),
\end{align*}
where $\tau=| t| +2.$ The implicit constants depend only on $%
\sigma.$
\end{lemma}

Set $\sigma=\frac{3}{2}$ and $x=T^{\alpha}.$ One computes the
integral of
the left-hand side:
\begin{equation*}
\int_{0}^{T}\left\vert2\sum_{\gamma}\frac{x^{i\gamma}}{1+\left(
t-\gamma
\right) ^{2}}\right\vert^{2}dt=2\pi F\left( \alpha,T\right) T\log
T+O\left( \left( \log T\right) ^{3}\right).
\end{equation*}
For the corresponding integral of the right-hand side, one uses the
Montgomery-Vaughan formula,
\begin{equation*}
\int_{0}^{T}\left\vert\sum_{n=1}^{\infty}\frac{a_{n}}{n^{it}}\right
\vert
^{2}dt=\sum_{n=1}^{\infty}\left\vert a_{n}\right\vert^{2}\left(
T+O\left(
n\right) \right)
\end{equation*}
and finds that the integral of the right-hand side equals
\begin{equation*}
T\log x+O\left( x\log x\right) +\frac{T}{x^{2}}\left[ \left( \log
T\right)
^{2}+O\left( \log T\right) \right],
\end{equation*}
which gives the claim of the theorem when one substitutes $x=T^{\alpha}.$

For more information about Montgomery conjecture, see \cite{goldston05}.

Montgomery's result allows us to compute the statistic (\ref%
{pair_correlations}) for pairs of Riemann's zeta zeros, provided that the
Fourier transform of the test function $r$ is supported on the interval
$%
[ -1,1].$ What about statistics of other $\ k$-tuples of the
zeros? First of all, the case of linear statistics ($k=1$) is similar
to the
questions considered by Selberg and Bourgade except that now we allow for
more general test functions. The main interest here is to see how far
we can
go in localizing this functions.

In this directions Hughes and Rudnick in \cite{hughes_rudnick02}
studied the
distribution of
\begin{equation*}
\mathcal{N}_{f}\left( t,T\right) :=\sum_{\gamma_{j}}f\left( \frac{\log
T}{%
2\pi}\left( \gamma_{j}-t\right) \right),
\end{equation*}
where $\gamma_{j}=\frac{1}{i}( \rho_{i}-\frac{1}{2}) $ and $%
\gamma_{j}$ are not assumed real. The function $f$ is a real-valued even
function with the smooth compactly-supported Fourier transform. (If $f$ is
the indicator function of an interval $[ -a,a] $ and if all $%
\gamma_{j}$ are real, then $\mathcal{N}_{f}( t,T) $ counts
number of zeros in the interval $[ t-a\frac{2\pi}{\log T},t+a\frac{%
2\pi}{\log T}].$ However, the Fourier transform of the indicator
function does not have compact support.)

Choose a weight function $w( x),$ such that $w\geq0,$ $\int
w( x) dx=1,$ and $\widehat{w}( x) $ is compactly
supported, and define an averaging operator
\begin{equation*}
\left\langle F\right\rangle_{T,H}:=\int_{\mathbb{R}}F\left( t\right)
w\left( \frac{t-T}{H}\right) \frac{dt}{H}.
\end{equation*}

\begin{theorem}[Hughes-Rudnick]
\label{theorem_Hughes_Rudnick1} Let the averaging window $H=T^{a}$ for $
0<a\leq1,$ and let be such that $\widehat{f}( u) =\int f(
x) e^{-2\pi ixu}dx\in C_{c}^{\infty}( \mathbb{R}) $ and $%
\mathrm{Supp}$ $f\subset( -2a/m,2a/m) $ with integer $m\geq1.$
Then, as $T\rightarrow\infty,$ the first $m$ moments of $\mathcal{N}_{f},$
$\langle \mathcal{N}_{f}^{m}\rangle_{T,H}$ converge to those of a Gaussian
random variable with expectation $\int f( x) dx$ and variance
\begin{equation*}
\sigma_{f}^{2}=\int\min\left( \left\vert u\right\vert,1\right)
\widehat{f%
}\left( u\right) ^{2}du.
\end{equation*}
\end{theorem}

Hence, if the frequency of the test function oscillations is bounded (and
therefore the function is very smooth and well delocalized in the $x$%
-space), then the first moments of the linear statistic converge to those
of the Gaussian variable. What about higher moments? Hughes and Rudnick show
(Theorem 6.5 in \cite{hughes_rudnick03}) that a similar result holds in the
random matrix theory for eigenvalues of a unitary random matrix. For random
matrices, higher moments do not converge to Gaussian values (Theorem
7.4 in
\cite{hughes_rudnick03}) Based on this analogy, they conjecture that it is
the same for the linear statistics of the zeta function zeros.

We will describe the ideas of the proof of Theorem \ref%
{theorem_Hughes_Rudnick1} below in the case when they are applied to the zeros
of Dirichlet's L-functions.

What about statistics of $k$-tuples of zeros when $k>2$? This case was
considered by Rudnick and Sarnak in \cite{rudnick_sarnak96}. Their results
hold for a quite large class of $L-$functions, and we will postpone their
discussion to a later section. Briefly, they are similar to Montgomery's
results since they show that the behavior of the zeta zeros is very similar
to the behavior of eigenvalues of random unitary matrices. Another similarity is that 
the results are
proven under some restrictive conditions on the Fourier transform of the
test function. It is an outstanding problem to prove that all results about
correlations of zeros hold without these restrictive hypotheses.

%s2.2 ###
\subsection{Dirichlet's L-functions}

In order to understand the behavior of the Riemann zeta zeros, it is
worthwhile to check for which functions their zeros have similar behavior. The simplest
example of such a family of functions is provided by Dirichlet's L-functions.

Let $\chi( n) $ denote a multiplicative \emph{character} modulo
a positive integer $q.$ That is, the function $\chi$ maps integers to the
unit circle; it is multiplicative, $\chi( nm) =\chi(
n) \chi( m),$ and $\chi( n) =0$ if $n$ and $q$
are not relatively prime. The character which sends every integer relatively
prime to $q$ to $1$ is called the \emph{principal character} modulo
$q$. The
\emph{conductor} of the character is the minimal integer $N$ such that the
character is periodic modulo $N.$ For simplicity, let $q$ be a prime in the
following. In this case the conductor equals $q.$ A character is \emph{odd}
if $\chi( -1) =-1,$ and \emph{even} if $\chi( -1)
=1. $

\emph{The Dirichlet }$L$\emph{-function} corresponding to the character
$%
\chi$ is defined by the series
\begin{equation*}
L\left( s,\chi\right) =\sum_{n=1}^{\infty}\frac{\chi\left( n\right)}{n^{s}}.
\label{def_Dirichlet_L}
\end{equation*}
This function has the Euler product representation:
%
%e14 ###
\begin{equation}
\zeta\left( s\right) =\prod_{p}\left( 1-\frac{\chi\left( p\right)
}{p^{s}}%
\right) ^{-1}
\end{equation}
because of the multiplicativity of $\chi( n) $. Moreover, the
argument behind the relation to theta functions (\ref{theta_representation})
can be repeated and as a consequence, one finds that series (\ref{def_Dirichlet_L})can be
continued to a function which is meromorphic in the whole complex plane and
satisfies a functional equation. Namely, let $\mu=0$ if $\chi$ is
even and
$\mu=1$ if $\chi$ is odd. Define
\begin{equation*}
\Phi\left( s,\chi\right) =q^{\frac{1}{2}\left( s+\mu\right) }\pi
^{-\frac{%
1}{2}\left( s+\mu\right) }\Gamma\left( \frac{1}{2}\left( s+\mu\right)
\right) L\left( s,\chi\right),
\end{equation*}
then the functional equation has the form
\begin{equation*}
\Phi\left( 1-s,\overline{\chi}\right) =\frac{i^{\mu}\sqrt{q}}{\tau
\left(
\chi\right) }\Phi\left( s,\chi\right),
\end{equation*}
where $\tau\left( \chi\right) $ is the Gauss sum:
\begin{equation*}
\tau\left( \chi\right) =\sum_{m=1}^{q}\chi\left( m\right) e^{2\pi im/q}.
\end{equation*}
The reason for appearance of $\tau(\chi)$ is that the modularity relation
(\ref%
{modularity_condition}) becomes more complicated in this case. For proofs,
see Chapter 9 in Davenport \cite{davenport67}.

Many other properties of the Dirichlet L-functions is similar to that
of the
Riemann zeta functions. In particular, one can establish similar explicit
formulas.

Two notable differences from the Riemann zeta is that (i) if $\chi$ is not
principal, then $L( s,\chi) $ is entire; and (ii) if $\chi$ is
not even, then a complex conjugate of a zero is not necessarily a zero.

%s2.2.1 ###
\subsubsection{Global scale}

For $T>0,$ let $N( T,\chi) $ denote the number of zeros of $
L( s,\chi) $ with $0<\sigma<1$ and $0\leq t\leq T,$ counting possible
zeros with $t=0$ or $t=T$ as one half only. Let
\begin{equation*}
S\left( t,\chi\right) =\frac{1}{\pi}\mathrm{Im}\log L\left(\frac
{1}{2}+it,\chi\right).
\end{equation*}
Then it can be shown that
\begin{equation*}
\mathcal{N}(T,\chi)=\frac{T}{2\pi}\log\frac{Tq}{2\pi e}-\frac{\chi
\left(
-1\right) }{8}+S\left( T,\chi\right) -S\left( 0,\chi\right) +O\left(
\frac{%
1}{1+T}\right).
\end{equation*}
See formula (1.3) in Selberg's paper \cite{selberg46} and Chapter 16 in
Davenport \cite{davenport67}.

If the character is fixed and $T$ is large, then\ the results for
$\mathcal{N%
}(T,\chi)$ are quite similar to results for $\mathcal{N}(T).$

A different situation arises when the interval $[ 0,T] $ is fixed
and the character $\chi$ varies (in particular, if $\chi$ is random).
\
This situation was studied by Selberg, who proved the following result (cf.
Theorem 9 in \cite{selberg46}).

\begin{theorem}[Selberg]
\label{theorem_selberg_L}For $\vert t\vert\leq
q^{1/4-\varepsilon},$ we have
\begin{equation*}
\frac{1}{q-2}\sum_{\chi}\left\vert S\left( t,\chi\right) \right\vert
^{2r}=%
\frac{\left( 2r\right) !}{r!\left( 2\pi\right) ^{2r}}\left( \log\log
q\right) ^{r}+O(\left( \log\log q\right) ^{r-1}),
\end{equation*}
where the summation if over all non-principal characters over the base $q.$
\end{theorem}

In other words, if $t$ is fixed and $q$ grows to infinity, then the
distribution of $S( t,\chi) $ approaches the distribution of a
Gaussian random variable with the variance $\frac{1}{2\pi^{2}}\log
\log q.$

One is naturally let to the question of correlations between $S( t,\chi
) $ for different $\chi.$ One result in this direction is stated by
Fujii (see p. 233 in \cite{fujii99}). Namely,
\begin{equation*}
\int_{0}^{T}S\left( t,\chi_{1}\right) S\left( t,\chi_{2}\right)
dt=\frac{%
\delta_{\chi_{1},\chi_{2}}}{2\pi}T\log T+A\left( \chi_{1},\chi
_{2}\right) T+O\left( \frac{T}{\sqrt{\log T}}\right),
\end{equation*}
where $A\left( \chi_{1},\chi_{2}\right)$ is a constant that depends on $\chi_{1}$ and $\chi_{2},$ 
which basically says that $S(t,\chi_{1}) $ and $S(t,\chi_{2}) $ 
are uncorrelated as functions of a random $t$ if $\chi_{1}\neq\chi_{2}.$ (See, however, a critique of Fujii's proof on p. 4 in \cite
{korolev10b}.)

Apparently, the question of correlations between $S( t_{1},\chi)
$ and $S( t_{2},\chi) $ as functions of a random $\chi$ has not yet been
investigated.

%s2.2.2 ###
\subsubsection{Local scale}

In \cite{hughes_rudnick03}, Hughes and Rudnick study the linear statistics
of low-lying zeros of $L$-functions on the local scale. Hughes and Rudnick
order the zeros $\rho_{i,\chi}=\frac{1}{2}+i\gamma_{i,\chi}$ as follows:
\begin{equation*}
\cdots\leq\mathrm{Re}\gamma_{-2,\chi}\leq\mathrm{Re}\gamma
_{-1,\chi
}<0\leq\mathrm{Re}\gamma_{1,\chi}\leq\mathrm{Re}\gamma_{2,\chi
}\leq
\cdots
\end{equation*}
and define
\begin{equation*}
x_{i,\chi}=\frac{\log q}{2\pi}\gamma_{i,\chi}.
\end{equation*}
Then they define
\begin{equation*}
W_{f}\left( \chi\right) =\sum_{i=-\infty}^{\infty}f\left( x_{i,\chi
}\right)
\end{equation*}
where $f$ is a rapidly decaying test function.

The question is to understand the behavior of the averages
\begin{equation*}
\left\langle W_{f}^{m}\right\rangle=\frac{1}{q-2}\sum_{\chi\neq\chi
_{0}}W_{f}\left( \chi\right).
\end{equation*}

The basis for their analysis is a variant of the explicit formula (\ref%
{formula_delsarte}) that relates a sum over zeros of $L( s,\chi)
$ to a sum over prime powers. This formula is a particular version of the
formula from Rudnick and Sarnak \cite{rudnick_sarnak96}, which is valid for
a more general class of zeta functions. Let $h( r) $ be any even
analytic function in the strip $-c\leq\mathrm{Im}r\leq1+c$ (for $c>0$)
such that $\vert h( r) \vert\leq A( 1+\vert
r\vert) ^{-( 1+\delta) }$ (for $r\in\mathbb{R},$ $%
A>0,$ $\delta>0$). Then (cf. formula (2.1) in \cite{hughes_rudnick03}),
%
%e15 ###
\begin{eqnarray}
\sum_{j}h\left( \gamma_{j,x}\right) &=&\frac{1}{2\pi}\int_{-\infty
}^{\infty}h\left( r\right) \left( \log q+G_{\chi}\left( r\right)
\right) dr
\label{formula_Rudnick_Sarnak0} \\
&&{}-\sum_{n}\frac{\Lambda\left( n\right) }{\sqrt{n}}\widehat{h}\left(
\log
n\right) \left( \chi\left( n\right) +\overline{\chi}\left( n\right)
\right), \notag
\end{eqnarray}
where $\widehat{h}( u) =\frac{1}{2\pi}\int h(r)e^{-iru}dr,$ and
\begin{equation*}
G_{\chi}\left( r\right) =\frac{\Gamma^{\prime}}{\Gamma}\left( \frac
{1}{2}%
+\mu\left( \chi\right) +ir\right) +\frac{\Gamma^{\prime}}{\Gamma
}\left(
\frac{1}{2}+\mu\left( \chi\right) -ir\right) -\frac{1}{2}\log\pi.
\end{equation*}
(Recall that $\mu( \chi) =0,$ if $\chi$ is even, and
$=1,$ if $%
\chi$ is odd.)

Take $h( r) =f( \frac{\log q}{2\pi}r) $, so that $%
\widehat{h}( u) =\frac{1}{\log q}\widehat{f}( \frac{u}{\log q%
}).$ We say that $f$ is \emph{admissible}, if it is a real, even
function whose Fourier transform $\widehat{f}( u) :=\int
f(r)e^{-2\pi iru}dr$ is compactly supported, and such that $\vert
f( r)\vert\leq A( 1+\vert r\vert)
^{-( 1+\delta) }.$ Then, from (\ref{formula_Rudnick_Sarnak0}) we
get the following decomposition:
\begin{equation*}
W_{f}\left( \chi\right) =\overline{W_{f}}\left( \chi\right)
+W_{f}^{osc}\left( \chi\right),
\end{equation*}
where
\begin{equation*}
\overline{W_{f}}\left( \chi\right) :=\int_{-\infty}^{\infty}f\left(
\frac{%
\log q}{2\pi}r\right) \left( \log q+G_{\chi}\left( r\right) \right) dr,
\end{equation*}
and
\begin{equation*}
W_{f}^{osc}\left( \chi\right) :=-\frac{1}{\log q}\sum_{n}\frac{\Lambda
\left( n\right) }{\sqrt{n}}\widehat{f}\left( \frac{\log n}{\log q}\right)
\left( \chi\left( n\right) +\overline{\chi}\left( n\right) \right).
\end{equation*}

For large $q,$
\begin{equation*}
\overline{W_{f}}\left( \chi\right) :=\int_{-\infty}^{\infty}f\left(
x\right) dx+O\left( \frac{1}{\log q}\right),
\end{equation*}
which is asymptotically independent of $\chi.$

For the oscillating part one has the following result (cf. Theorem 5.1\ in
\cite{hughes_rudnick03}).

\begin{theorem}[Hughes and Rudnick]
Let $f$ be an admissible function and assume that
\begin{equation*}
\mathrm{Supp}\left( \widehat{f}\right) \subseteq\left[ -\alpha,\alpha
\right],\text{ }\alpha>0.
\end{equation*}
If $m<2/\alpha,$ then the $m$-th moment of $W_{f}^{osc}$ is \
\begin{equation*}
\lim_{q\rightarrow\infty}\left\langle\left( W_{f}^{osc}\right)
^{m}\right\rangle_{q}=\left\{
\begin{array}{cc}
\frac{m!}{2^{m/2}\left( m/2\right) !}\sigma\left( f\right) ^{m}, &
\text{if
}m\text{ is even,} \\
0, & \text{if }m\text{ is odd,}%
\end{array}
\right.
\end{equation*}
where
\begin{equation*}
\sigma\left( f\right) ^{2}=\int_{-1}^{1}\left\vert u\right\vert
\widehat{f}%
\left( u\right) ^{2}du.
\end{equation*}
\end{theorem}

In other words, the first several moments of the statistic $W_{f}$ converge
to the corresponding moments of a Gaussian random variable. A similar
situation holds for an eigenvalue statistic of random unitary matrices.
Hughes and Rudnick show (Theorem 7.4) that higher moments for this
eigenvalue statistic are not Gaussian (using results from the work of
Diaconis and Shahshahani \cite{diaconis_shahshahani94}). They conjecture
that the same result should hold for the statistic $W_{f}.$

As an application, Hughes and Rudnick derived some results for the smallest
zero of $L( s,\chi).$ In particular, they showed that for
infinitely many $q$ there are characters $\chi$ such that the imaginary
part of the zero is between $0$ and $1/4$ (Corollary 8.2 in \cite%
{hughes_rudnick03}). They conjecture that $1/4$ can be substituted with
arbitrary positive constant.\vadjust{\goodbreak}

Moreover, if $\beta>6.333,$ then a proportion of characters $\chi$
whose $%
L $ function has a zero with imaginary part between $\ 0$ and $\beta$ is
greater than $c( \beta) >0$ for all sufficiently large $q.$
(Theorem 8.3 in \cite{hughes_rudnick03}). The conjecture is that this
is in
fact true for every $\beta>0$.

%s2.3 ###
\subsection{L-functions for modular forms (The Hecke L-functions)}

Hecke generalized the Riemann zeta function by using ideals in an imaginary
quadratic field $K$ instead of integers:
\begin{equation*}
L_{K}\left( s\right) =\sum_{\mathfrak{a}}\left( N\left( \mathfrak
{a}\right)
\right) ^{-s}.
\end{equation*}
This function has a Euler product formula and a functional equation,
although the latter is somewhat different: let
\begin{equation*}
\Lambda_{K}\left( s\right) =\left( \frac{\sqrt{\left\vert D\right\vert
}}{%
2\pi}\right) ^{s}\Gamma\left( s\right) L_{K}\left( s\right),
\end{equation*}
where $D$ is the discriminant of the imaginary quadratic field $K.$ Then
\begin{equation*}
\Lambda_{K}\left( 1-s\right) =\Lambda_{K}\left( s\right).
\end{equation*}
(Compare this with (\ref{functional_equation_Riemann}).) The proof of this
functional equation is similar to Riemann's proof and relies on a modularity
property of a certain complex-analytic function, that is, on its behavior
relative to the change of variable $z\rightarrow1/z.$

Motivated by this example, Hecke had a fruitful idea of obtaining
L-functions from complex-analytic functions that transform well under the
action of the modular group, and then checking which additional conditions
are needed to ensure that a functional equation and a Euler product formula
holds. The objects constructed in this way are called Hecke L-functions.

In order to illustrate, let
\begin{eqnarray*}
f\left( z\right) &=&\sum_{n>0}c_{n}e^{2\pi inz},\text{ and } \\
L\left( f,s\right) &=&\sum_{n>0}c_{n}n^{-s},
\end{eqnarray*}
Define
\begin{equation*}
\Lambda\left( f,s\right) =\left( \frac{\sqrt{N}}{2\pi}\right)
^{s}\Gamma
\left( s\right) L\left( f,s\right).
\end{equation*}
Since
\begin{equation*}
\left( \frac{\sqrt{N}}{2\pi}\right) ^{s}\Gamma\left( s\right)
n^{-s}=\int_{0}^{\infty}e^{-2\pi ny/\sqrt{N}}y^{s}\frac{dy}{y},
\end{equation*}
we have the representation
\begin{equation*}
\Lambda\left( f,s\right) =\int_{0}^{\infty} f\left( \frac
{iy}{\sqrt{N%
}}\right) y^{s}\frac{dy}{y},\vadjust{\goodbreak}
\end{equation*}
which implies that
%
%e16 ###
\begin{equation}
\Lambda\left( f,s\right) =\int_{1}^{\infty} f\left( \frac
{iy}{\sqrt{N%
}}\right) y^{s}\frac{dy}{y}+i^{k}\int_{1}^{\infty}Wf\left( \frac{iy}{\sqrt{N}}\right)  y^{k-s}\frac{dy}{y},
\label{Hecke_representation}
\end{equation}
where
%
%e17 ###
\begin{equation}
Wf\left( z\right) =\left( \sqrt{N}z\right) ^{-k}f\left( \frac
{-1}{Nz}\right)
. \label{Fricke_involution}
\end{equation}
Consequently, any eigenfunction of the operator $W$ will have a functional
equation similar to the functional equation for the Riemann zeta function.
The idea is to find a suitable finite-dimensional space of functions $f,$
which is invariant under the action of $W$ and diagonalize $W$ in this
finite-dimensional space.

Below, we give a brief outline of these ideas. A good source for this
material is Chapter 14 in \cite{iwaniec_kowalski04} and Chapter V in
\cite%
{milne06}.

%s2.3.1 ###
\subsubsection{L-functions from modular forms}

\label{section_L_modular_forms}

Let $\Gamma$ is a subgroup of finite index in $SL_{2}( \mathbb{Z}%
) $ and let $\mathbb{H}$ denote the upper half-plane $\{ z|%
\mathrm{Im}z>0\}.$ The set $\mathbb{H}^{\ast}=\mathbb{H\cup}\{
\infty\} \cup\mathbb{Q}$ can be made into a Hausdorff topological
space and one can define a continuous action of $SL_{2}( \mathbb{Z}%
) $ on $\mathbb{H}^{\ast}$ as an extension of the action of $%
SL_{2}( \mathbb{Z}) $ on $\mathbb{H}:$
\begin{equation*}
\text{if }\gamma=\left(
\begin{array}{cc}
a & b \\
c & d%
\end{array}
\right),\text{ then }\gamma z=\frac{az+b}{cz+d}.
\end{equation*}
The \emph{cusps} are points in $\mathbb{H}^{\ast}\backslash\mathbb{H}$.
One can also show that $\Gamma\backslash\mathbb{H}^{\ast}$ is a compact
Hausdorff space (that is, a compact space in which every two points have
disjoint open neighbourhoods), which is a Riemann surface (that is, it
admits a complex structure).

Next, we define an action of $SL_{2}( \mathbb{Z}) $ on functions $%
f:\mathbb{H}^{\ast}\rightarrow\mathbb{C}:$
\begin{equation*}
\text{if }\gamma=\left(
\begin{array}{cc}
a & b \\
c & d%
\end{array}
\right),\text{ then }f\circ\left[ \gamma\right] _{k}:=\left(
cz+d\right)
^{-k}f\left( \gamma z\right).
\end{equation*}
In fact if $\gamma\in GL_{2}( \mathbb{Z}),$ then one can define
its action by
\begin{equation*}
f\circ\left[ \gamma\right] _{k}:=f\circ\left[ \det\left( \gamma
\right)
^{-1/2}\gamma\right] _{k}.
\end{equation*}
It can be checked that this is indeed a group action of $GL_{2}(
\mathbb{Z}) $ on functions. 

While modular forms can be defined for any discrete subgroup $\Gamma,$ the
most studied are subgroups $\Gamma_{0}( N) $:
\begin{equation*}
\Gamma_{0}\left( N\right) :=\left\{ \left(
\begin{array}{cc}
a & b \\
c & d%
\end{array}
\right) |\text{ }c\equiv0\left( \mathrm{mod}N\right) \right\}.
\end{equation*}
For these subgroups, we have the following definition.

\begin{definition}
Let $\chi$ be a Dirichlet character modulo $N.$ A \emph{modular form}
for $%
\Gamma_{0}( N) $ of weight $k\geq1$ and character $\chi$ is a
function $\mathbb{H}^{\ast}\rightarrow\mathbb{C}$ such that
\begin{itemize}
\item[(a)] $f$ is holomorphic on $\mathbb{H}$;
\item[(b)] for any $\gamma\in\Gamma_{0}( N),$ $f\circ[ \gamma%
] _{k}=\chi( d) f;$
\item[(c)] $f$ is holomorphic at the cusps.
\end{itemize}
A \emph{cusp form} is a modular form which is zero at all the cusps.
\end{definition}

In particular, the complex analyticity at infinity implies that a modular
form for $\Gamma_{0}( N) $ can be written as
\begin{equation*}
f\left( z\right) =\sum_{n\geq0}c_{n}e^{2\pi inz}.
\end{equation*}
For a cusp form, $c_{0}=0.$

For simplicity, \textit{we will only consider the forms with the principal character
$\chi$} and we will write $\mathcal{M}_{k}( N) $ and $\mathcal{S}%
_{k}( N) $ to denote the linear spaces of the modular and cusp
forms. One can show that these spaces are finite dimensional for each~$k.$
Note that if $k$ is odd then $\mathcal{M}_{k}( N) $ is zero since
$f\circ[ -I] _{k}=( -1) ^{k}f,$ hence we should have $%
f=-f.$

\begin{definition}
Let $f$ be a cusp form of weight $2k$ for $\Gamma_{0}( N),$
\begin{equation*}
f\left( z\right) =\sum_{n>0}c_{n}e^{2\pi inz}.
\end{equation*}
The $L$-series of the cusp form $f$ \ is the Dirichlet series
\begin{equation*}
L\left( f,s\right) =\sum_{n>0}c_{n}n^{-s}.
\end{equation*}
\end{definition}

It is possible to estimate that $\vert c_{n}\vert\leq Cn^{k},$
and therefore this series is convergent for $\mathrm{Re}s>k+1.$

The crucial fact is that the space of modular forms is invariant under the
action of operator $W$ from (\ref{Fricke_involution}). Indeed,
$Wf=f\circ%
[ \omega] _{k},$ where $\omega=\left(
\begin{smallmatrix}
0 & -1 \\
N & 0%
\end{smallmatrix}
\right).$ If $\gamma=\left(
\begin{smallmatrix}
a & b \\
c & d%
\end{smallmatrix}%
\right) \in\Gamma_{0}( N).$ Then
\begin{equation*}
\omega\gamma\omega^{-1}=\left(
\begin{array}{cc}
d & -c/N \\
-Nb & a%
\end{array}
\right) \in\Gamma_{0}\left( N\right).
\end{equation*}
Hence
\begin{equation*}
Wf\circ\left[ \gamma\right] _{k}=f\circ\left[ \omega\gamma\omega
^{-1}%
\right] _{k}\circ\left[ \omega\right] _{k}=f\circ\left[ \omega\right]
_{k}=Wf.
\end{equation*}
where the second step is by modularity of $f.$ One can also check that
\thinspace$W$ preserves $\mathcal{S}_{2k}( N) $. Since $W^{2}=1,$
hence the only eigenvalues of $w_{N}$ are $\pm1,$ and $\mathcal{S}%
_{2k}( N) $ is a direct sum of the corresponding eigenspaces, $%
\mathcal{S}_{2k}=\mathcal{S}_{2k}^{+1}+\mathcal{S}_{2k}^{-1}$.

By the argument in the beginning of this section (see formula (\ref%
{Hecke_representation})), one can infer the following result.

\begin{theorem}[Hecke]
Let $f\in\mathcal{S}_{2k}( \Gamma_{0}( N) ) $ be a
cusp form in the $\varepsilon$-eigenspace, $\varepsilon=1$ or $-1.$ Then
the function $\Lambda( s,f) :=N^{s/2}( 2\pi)
^{-s}\Gamma( s) L( f,s) $ extends analytically to a
holomorphic function on the whole complex plane, and satisfies the
functional equation
\begin{equation*}
\Lambda\left( s,f\right) =\varepsilon\left( -1\right) ^{k}\Lambda
\left(
k-s,f\right).
\end{equation*}
\end{theorem}

The natural question is what about the Euler product formula?

By studying the properties of the modular forms that arise from the $L$%
-functions of the quadratic imaginary fields (and that have the product
formula almost by definition), Hecke was able to formulate a list of
properties which should be imposed on the modular form $f,$ so that
$L(f,z) $ had a product formula.

Namely, define the Hecke operators (cf. formula (14.46) in \cite%
{iwaniec_kowalski04})
\begin{equation*}
\left[ T\left( n\right) f\right] \left( z\right) :=\frac{1}{n}%
\sum_{ad=n}a^{k}\sum_{0\leq b<d}f\left( \frac{az+b}{d}\right).
\end{equation*}
It can be checked that these are linear operators on $\mathcal
{S}_{2k}(\Gamma_{0}( N) ).$ (See Section IX.6 in Knapp \cite%
{knapp92} or Section V.4 of Milne \cite{milne06} or Section VII.5 of Serre
\cite{serre73} for details). They have the following properties:

\begin{theorem}[Hecke]
The maps $T( n) $ have the following properties:
\begin{itemize}
\item[(a)] $T( mn) =T( m) T( n) $ if $\gcd (m,n) =1;$
\item[(b)] $T( p) T( p^{r}) =T( p^{r+1})
+p^{2k-1}T( p^{r-1}) $ if $p$ does not divide $N;$
\item[(c)] $T( p^{r}) =T( p) ^{r},$ $r\geq1,$ if $p|N;$%
\item[(d)] all $T( n) $ commute.
\end{itemize}
\end{theorem}

Moreover, by acting on the Fourier expansion, one finds that the first
Fourier coefficient in the expansion of $T( n) f$ is $c_{n}$.
Hence, if $f$ is an eigenfunction of $f$ with eigenvalue $\lambda_{n},$
then $c_{n}=\lambda_{n}c_{1}.$

Since the Hecke operators commute we can look for the modular functions $f$
which are eigenfunctions for all of them. Then, the multiplicativity properties
of $T( n) $ imply the corresponding properties for
coefficients $%
c_{n},$ which leads to a Euler product formula for $L( f,s) $.
This is formalized in the following result.

\begin{theorem}[Hecke]
Let $f$ be a cusp form of weight $2k$ for $\Gamma_{0}( N) $ that
is simultaneously an eigenvector for all $T( n),$ say
$T(n) f=\lambda_{n}f,$ and let
\begin{equation*}
f=\sum_{n\geq1}c_{n}q^{n},\text{ }q=e^{2\pi iz}.
\end{equation*}
Let $c_{1}=1.$ Then, (i) coefficients of the series are eigenvalues of the
Hecke operators,
\begin{equation*}
c_{n}=\lambda_{n}
\end{equation*}
and (ii)
\begin{equation*}
L\left( f,s\right) =\prod_{p|N}\frac{1}{1-c_{p}p^{-s}}\prod_{p\nmid
N}\frac{1%
}{1-c_{p}p^{-s}+p^{2k-1-s}}.
\end{equation*}
\end{theorem}

For example, $\mathcal{S}_{12}( \Gamma_{0}( 1) ) $
has dimension $1,$ and therefore it is generated by a single function, which
is called the $\Delta$-function:
\begin{equation*}
\Delta\left( q\right) =q\prod_{1}^{\infty}\left( 1-q^{n}\right)
^{24}=\sum
\tau\left( n\right) q^{n},
\end{equation*}
where $\tau( n) $ is the Ramanujan $\tau$-function. It follows
that the $L$-function associated to the $\Delta$-function has both a
functional equation and the Euler product property.

In a more general situation, if we wish to find forms that have both a
functional equation and a Euler product, then we must overcome the obstacle
that in some exceptional cases operators $W$ and $T( n) $ do not
commute. However, this obstacle can be circumvented and it can be proved
that such good modular forms do exist. They are called \emph{primitive
forms} or
\emph{newforms. }

In summary, the L-functions of primitive forms have both a functional
equation and the Euler product property. As a consequence, one can write
explicit formulas for these L-functions.

%s2.3.2 ###
\subsubsection{L-functions from Maass forms}

A nice source for the material in this section is the lecture notes by Liu
\cite{liu07}. A~lot of additional information about Maass forms can be found
in the book by Iwaniec \cite{iwaniec02}.

Modular forms are holomorphic and they are not easy to construct or compute.
One can try to use Hecke ideas for a different class of functions that
satisfy a modularity condition. In this way one comes to the concept of a
Maass form.

\begin{definition}
A smooth function $f\neq0$ is called a \emph{Maass form} for group
$\Gamma
, $ if
\begin{itemize}
\item[(i)] for all $g\in\Gamma$ and all $z\in\mathbb{H},$ $f( gz)
=f( z);$\newline
\item[(ii)] $f$ is an eigenfunction of the non-Euclidean Laplace operator:
\begin{equation*}
-y^{2}\left( \frac{\partial^{2}}{\partial x^{2}}+\frac{\partial^{2}}{%
\partial y^{2}}\right) f=\lambda f,
\end{equation*}
and
\item[(iii)] there exists a positive integer $N,$ such that
\begin{equation*}
f\left( z\right) \ll y^{N},\text{ }y\rightarrow\infty.
\end{equation*}
\end{itemize}
A Maass form f is said to be a \emph{cusp form} if the equality
\begin{equation*}
\int_{0}^{1}f\left( z+b\right) db=0
\end{equation*}
holds for all $z\in\mathbb{H}$.
\end{definition}

A Maass form $f$ is call \emph{odd} if $f( -x+iy) =-f(
x+iy),$ and \emph{even} if $f( -x+iy) =f(
x+iy).$\vadjust{\goodbreak}

Note that it is relatively easy to generate Maass forms as
eigenfunctions of
the Laplace operator on a fundamental domain of the group $\Gamma.$ By
expanding a Maass form in Fourier series and taking the Fourier coefficients
as the coefficients of a Dirichlet series, one can construct new
L-functions. Precisely, let $f$ be either an even or an odd cusp Maass form
with eigenvalue $1/4+r^{2}$. Then, one can write:\looseness=-1
\begin{equation}
f\left( x+iy\right) =\sqrt{y}\sum_{n\neq0}c_{n}K_{ir}\left( 2\pi\left
\vert
n\right\vert y\right) e^{2\pi inx},
\label{Maass_form_expansion}
\end{equation}\looseness=0
where $K_{ir}$ are Bessel functions, and we define
%
%e18 ###
\begin{equation}
L\left( f,s\right) =\sum\limits_{n>0}c_{n}n^{-s}. \label{Maass_L}
\end{equation}
The key idea here is the fact that the Laplace operator commutes with Hecke
operators, and therefore all these operators can be simultaneously
diagonalized. By a computation, the first Fourier coefficient of $T(
n) f$ is $c_{n}c_{1}.$ As a consequence, $L$-functions corresponding
to Maass forms have a product formula.

What about the functional equation? It holds. However, instead of the
standard formula for the Gamma function one needs the following integral:
\begin{equation*}
\int_{0}^{\infty}K_{ir}\left( y\right) y^{s}\frac{dy}{y}=\Gamma\left(
\frac{s+ir}{2}\right) \Gamma\left( \frac{s-ir}{2}\right).
\end{equation*}

\begin{theorem}
Let $f$ be a Maass form with eigenvalue $1/4+r^{2}.$ Let $\varepsilon
=0$ or
$1$ depending on whether $f$ is even or odd. Let
\begin{equation*}
\Lambda\left( f,s\right) =\pi^{-s}\Gamma\left( \frac{s+\varepsilon
+ir}{2}%
\right) \Gamma\left( \frac{s+\varepsilon-ir}{2}\right) L\left(
f,s\right).
\end{equation*}
Then $\Lambda( s,f) $ is an entire function that satisfies
\begin{equation*}
\Lambda\left( f,s\right) =\left( -1\right) ^{\varepsilon}\Lambda\left(
f,1-s\right).
\end{equation*}
\end{theorem}

%s2.3.3 ###
\subsubsection{Statistical properties of zeros: Global scale}

Let
\begin{equation*}
S\left( f,t\right) :=\frac{1}{\pi}\arg L\left( f,\frac{1}{2}+it\right),
\end{equation*}
where $f$ is the Maass form with an eigenvalue $\lambda$ and $L$ is the
corresponding $L$-function. $S( f,t) $ is related to the number
of zeros of $L$ in the critical strip in the same way as the usual
$S(t) $ function is related to the number of zeros of Riemann's zeta
function.

We are interested here in the distribution of $S( f,t) $ with
respect to the random choice of $f$.

Of course one need to explain what is meant by the random choice of
$f.$ Let
$S_{j}( t) :=S( f_{j},t) $ where $f_{j}$ has an
eigenvalue $\lambda_{j}=\frac{1}{4}+r_{j}^{2}.$ Define $\nu_{j}(
n) :=c_{j}( n) /\sqrt{\cosh\pi r_{j}},$ where
$c_{j}(
n) $ are coefficients in the expansion (\ref{Maass_form_expansion}) for the Maass form 
$f_{j}.$ The numbers $\nu_{j}(1)$ will be used as
weights in the limiting procedure. (We assume that $f_{j}$ are normalized to have a
unit norm as $L^{2}$-functions. Therefore,\vadjust{\goodbreak} $\nu_{j}(
1) $ are not necessarily equal to $1.$)  One knows that $\nu_{j}( 1) $
are $O( \sqrt{r_{j}}) $ and
\begin{equation*}
\frac{1}{T^{2}}\sum_{r_{j}\leq T}\left\vert\nu_{j}\left( 1\right)
\right\vert^{2}=\frac{1}{\pi^{2}}+O\left( \frac{\log T}{T}\right).
\end{equation*}
Since by the Weyl law, the number of $r_{j}$ below $T$ is proportionate
to $%
T^{2},$ these weights can be thought as having bounded magnitude and
not too
sparse. Here is one of the results about the randomness of $S_{j}(
t) $ (Theorem 3 in \cite{hejhal_luo97}).

\begin{theorem}[Hejhal-Luo]
Let $h>0$ and $t>0$ be fixed. Then we have
\begin{equation*}
\lim_{T\rightarrow\infty}\frac{1}{\left( 2HT\right) }\sum_{\left\vert
r_{j}-T\right\vert\leq H}\frac{\pi^{2}\left\vert\nu_{j}\left(
1\right)
\right\vert^{2}}{2}\frac{\left( S_{j}\left( t\right) \right)
^{n}}{\left(
\log\log T\right) ^{n/2}}=C_{n},
\end{equation*}
where $C_{n}$ are moments of the Gaussian distribution.
\end{theorem}

%s2.3.4 ###
\subsubsection{Local scale}

\label{section_Rudnick_Sarnak}

Rudnick and Sarnak \cite{rudnick_sarnak96} extended the results of
Montgomery to zeta functions that arise from modular and Maass
forms. In fact, they work in greater generality and study the zeta functions
that arise from the automorphic cuspidal representations of
$GL_{m}.$%
 The Hecke modular L-functions correspond to the case $m=2.$ Their main
tool is the following explicit formula, which we formulate for the case of
the Hecke L-functions. Let
\begin{eqnarray*}
L\left( s,f\right) &=&\prod_{p|N}\frac{1}{1-c\left( p\right) p^{-s}}%
\prod_{p\nmid N}\frac{1}{1-c\left( p\right) p^{-s}+p^{2k-1-s}} \\
&=&\prod L_{p}\left( s,f\right).
\end{eqnarray*}
where
\begin{equation*}
L_{p}\left( s,f\right) =\frac{1}{\left( 1-\alpha_{1}(p)p^{-s}\right)
\left(
1-\alpha_{2}(p)p^{-s}\right) },
\end{equation*}
with the convention that for $p|N$ one of $\alpha_{i}( p) $ is
zero. Let $a( p^{k}) =\alpha_{1}( p) ^{k}+\alpha
_{2}( p) ^{k},$ and define $b( n) =\Lambda(
n) a( n).$ Then
\begin{equation*}
\frac{L^{\prime}}{L}=-\sum_{n=1}^{\infty}\frac{b\left( n\right) }{n^{s}}.
\end{equation*}

\begin{theorem}[Rudnick and Sarnak]
Let $\widehat{h}\in C_{c}^{\infty}( \mathbb{R}) $ be a smooth
compactly supported function, and let $h( r) =\int\widehat
{h}( u) e^{iru}du$. Then
\begin{align}
\sum h( \gamma) &=\frac{\log Q}{2\pi}\int_{-\infty
}^{\infty
}h(r) dr \notag\\
&\quad {}+\frac{1}{2\pi}\int_{-\infty}^{\infty
}h( r) \Biggl( \sum_{j}\biggl[ \frac{\Gamma^{\prime}}{%
\Gamma}\biggl( \frac{1}{2}+\mu_{j}+ir\biggr) +\frac{\Gamma^{\prime}}{
\Gamma}\biggl( \frac{1}{2}+\mu_{j}-ir\biggr) \biggr] \Biggr) dr \notag
\\
&\quad {}-\sum_{n=1}^{\infty}\biggl( \frac{b( n) }{\sqrt{n}}\widehat
{h}%
( \log n) +\frac{\overline{b( n) }}{\sqrt{n}}\widehat{h%
}( -\log n)\biggr), \label{formula_Rudnick_Sarnak}
\end{align}
where $\mu_{j}$ are some parameters that depend on the form $f,$ and
$Q$ is
the conductor of the form.
\end{theorem}

By using this result and estimates on the size of coefficients $b(
n) $, Rudnick and Sarnak proved a generalization of the Montgomery
theorem. Their result is valid not only for the Riemann zeta function, but
also for Dirichlet $L$-functions, for Hecke modular $L$-functions and
for $L$%
-functions that correspond to automorphic cuspidal representations of $%
GL_{3}.$ We formulate it for Hecke modular $L$-functions.

Consider the class of smooth test functions $F(x_{1},\ldots,x_{n})$ that
satisfy the following conditions:

\begin{description}
\item[\textbf{TF 1}] $F(x_{1},\ldots,x_{n})$ is symmetric.
\item[\textbf{TF 2}] $F(x+t( 1,\ldots,1) )=F( x) $ for
all $%
t\in\mathbb{R}$.
\item[\textbf{TF 3}] $F( x) \rightarrow0$ rapidly as $|x|\rightarrow\infty$
in the hyperplane $\sum x_{j}=0.$
\end{description}

If $B_{N}$ is a set of $N$ numbers $x_{1},\ldots,x_{N},$ then the $n$-level
correlation sum is defined by
\begin{equation*}
R_{n}\left( B_{N},F\right) =\frac{n!}{N}\sum_{\substack{ S\subset B_{N}
\\ %
\left| S\right| =n}}F\left( S\right).
\end{equation*}

Define the $n$-level correlation density by
\begin{equation*}
W_{n}\left( x_{1},\ldots,x_{n}\right) =\det\left( K\left(
x_{i}-x_{j}\right) \right),\text{ }K\left( x\right) =\frac{\sin\pi
x}{\pi x%
}.
\end{equation*}

Then the following result holds (cf. Theorem 1.2 in \cite{rudnick_sarnak96}).

\begin{theorem}
Assume the Riemann hypothesis for the zeros of $L( s,f).$
Let $%
F( x_{1},\ldots,x_{N}) $ satisfy \textbf{TF 1, 2, 3} and in
addition assume that $\widehat{F}( \xi) $ is supported in $%
\sum_{j}| \xi_{j}| <1.$ Then,
\begin{equation*}
R_{n}\left( B_{N},F\right) \rightarrow\int F\left( x\right) W_{n}\left(
x\right) \delta\left( \frac{x_{1}+\cdots+x_{n}}{n}\right) dx_{1}\ldots
dx_{n}
\end{equation*}
as $N\rightarrow\infty.$
\end{theorem}

Rudnick and Sarnak mention that the result can probably be proven for
functions $F$ with the Fourier transform supported in $\sum_{j}|
\xi_{j}| <2$ by an improvement of their method, and conjecture that it
holds without any assumption on the support of $\widehat{F}( \xi
).$

%s2.4 ###
\subsection{ Elliptic curve zeta functions}

The main source for this section is the book \cite{milne06} \ by Milne.
Consider an elliptic curve
\begin{equation*}
E:Y^{2}Z=X^{3}+aXZ^{2}+bZ^{3},
\end{equation*}
where $a$ and $b$ are integer, and assume that $\vert\Delta
\vert=\vert4a^{3}+27b^{2}\vert$ cannot be made smaller
by a change of variable $X\rightarrow X/c^{2},\,Y\rightarrow Y/c^{3}.$ This
equation is called \emph{minimal}. The equation
\begin{equation*}
\overline{E}:Y^{2}Z=X^{3}+\overline{a}XZ^{2}+\overline{b}Z^{3},
\end{equation*}
with $\overline{a}$ and $\overline{b}$ the images of $a$ and $b$ in
$\mathbb{%
F}_{p}$ (the finite field with $p$ elements) is called \emph{the reduction
of }$E$\emph{\ modulo }$p.$ (It is assumed here that $p\neq2,3.$ In the
case when $p$ is $2$ or $3$, a somewhat different notion of the minimal
equation is needed.) Let $N_{p}$ is the number of solutions of this equation
in $\mathbb{F}_{p}.$

There are four possible cases:
\begin{itemize}
\item[(a)] \textbf{Good reduction}. $\overline{E}$ is an elliptic curve. (That is,
the determinant does not vanish and therefore the curve is smooth.) This
happens if $p\neq2$ and $p$ does not divide $\Delta.$

\item[(b)] \textbf{Cuspidal, or additive, reduction.} This is the case in
which the
reduced curve has a cusp. For $p\neq2,3$, this case occurs exactly
when $%
p|4a^{3}+27b^{2}$ and $p|-2ab.$

\item[(c)] \textbf{Nodal, or multiplicative, reduction.} The reduced curve has a
node. For $p\neq2,3$, it occurs exactly when $p|4a^{3}+27b^{2}$,
$p\nmid
-2ab.$
\begin{itemize}
\item[(c1)] \textbf{Split case.} The tangents at the node are rational
over $%
\mathbb{F}_{p}.$ This happens when $-2ab$ is a square in $\mathbb{F}_{p}$.

\item[(c2)] \textbf{Non-split case.} The tangents at the node are not
rational over $\mathbb{F}_{p}.$This occurs when $-2ab$ is not a square
in $%
\mathbb{F}_{p}$.
\end{itemize}
\end{itemize}

The names additive and multiplicative refer to the group of points on the
reduced curve, which in these cases is isomorphic either to $(
\mathbb{F%
}_{p},+),$ or $( \mathbb{F}_{p}^{\ast},\times).$

We define the $L$ function associated with the elliptic curve $E$ as
follows.
\begin{equation*}
L\left( s,E\right) :=\prod_{p}\frac{1}{L_{p}\left( p^{-s}\right) }.
\end{equation*}

Here, the local factors $L_{p}( T) $ are defined as follows:
\begin{equation*}
L_{p}\left( T\right) =\left\{
\begin{array}{cc}
1-a_{p}T+pT^{2}, & \text{if }p\text{ is good, with }a_{p}=p+1-N_{p}, \\
1-T & \text{if }E\text{ has split multiplicative reduction,} \\
1+T & \text{if }E\text{ has non-split multiplicative reduction,} \\
1 & \text{if }E\text{ has additive reduction.}%
\end{array}
\right.
\end{equation*}

Let $S$ be the (finite) set of primes with bad reduction. Then we can also
write
\begin{eqnarray*}
L\left( s,E\right) &:=&\prod_{p\in S}\frac{1}{L_{p}\left( p^{-s}\right)
}%
\prod_{p\notin S}\frac{1}{1+\left( N_{p}-p-1\right) p^{-s}+p^{1-2s}} \\
&=&\prod_{p\in S}\frac{1}{L_{p}\left( p^{-s}\right) }\prod_{p\notin
S}\frac{1%
}{\left( 1-\alpha_{p}p^{-s}\right) \left( 1-\beta_{p}p^{-s}\right) }.
\end{eqnarray*}

The Hasse-Weil conjecture says that $L( s,E) $ can be
analytically continued to a meromorphic function on the whole of
$\mathbb{C}$
and satisfies a functional equation. A recent work by Wiles and others
confirmed this conjecture by showing that all elliptic curves are
``modular'', in particular, their $L$%
-functions arise from modular forms. To a certain extent, this result
reduces the study of the elliptic $L$-functions to the study of the Hecke
modular $L$-functions.

It is known that the numbers $a_{p}$ do not exceed $2\sqrt{p}$ in absolute
value. For a fixed elliptic curve and different primes $p,$ these numbers
are believed to be distributed on the interval $[ -2\sqrt{p},2\sqrt{p}%
] $ according to the semicircle distribution but this is not proven.
In fact, this conjecture is related to the Birch-Swinnerton conjecture that
states that
\begin{equation*}
L\left( s,E\right) \sim C\left( s-1\right) ^{r}\text{ as }s\rightarrow1,
\end{equation*}
where $r$ is the rank of the group of rational points on $E,$ and $C$
is a
certain predicted constant. (see Chapter 10 in \cite{milne06} for more
information.)

It is possible to construct zeta functions for other nonsingular projective
varieties and the conjecture by Hasse and Weil states that these zeta
functions satisfy a functional equation (and the Riemann hypothesis).
However, apparently not much is known beyond the cases of projective spaces
and elliptic curves.

%s3 ###
\section{Selberg's zeta functions for compact and non-compact manifolds}

It is useful to keep in mind that we will now talk about a new type of zeta
functions, which is significantly different from number-theoretical zeta
functions. While there is an explicit formula, it relates Laplace eigenvalues and
geodesics, not zeta zeros and primes. The possibility of a relation between
these two types of zetas is only conjectural.

The main source for this section is Hejhal's book \cite{hejhal76}.

%s3.1 ###
\subsection{Selberg's zeta function and trace formula}

Let $M$ be a compact Riemann surface of genus $g\geq2.$ Then, $M$ can be
identified with a quotient space $\Gamma\backslash\mathbb{H},$ where $
\mathbb{H}$ is the upper half-plane and $\Gamma$ is a discrete
subgroup of $%
SL_{2}( R) $. We assume that $\mathbb{H}$ has the Poincare metric
$\vert dz\vert/y$ with the area element $dxdy/y^{2},$ and
therefore it has the constant negative curvature. This metric is naturally
projected on the surface $M.$

This is not the most general situation of interest since most of the
quotient spaces $\Gamma\backslash\mathbb{H}$ occurring in arithmetic
applications have cusps and therefore are non-compact. However, the theory
is most clear and transparent for the compact surfaces.

The Laplace operator on $M$ can be defined by the following formula.
\begin{equation*}
-\Delta:u\rightarrow-y^{2}\left( u_{xx}+u_{yy}\right).
\end{equation*}
It can be shown that this operator has a discrete set of non-positive
eigenvalues:
\begin{equation*}
0=\lambda_{0}<\lambda_{1}\leq\lambda_{2}\leq\ldots,
\end{equation*}
and the only point of accumulation of these eigenvalues is
$\infty.$\vadjust{\goodbreak}

Let us define
\begin{equation*}
r_{n}=\left\{
\begin{array}{cc}
\sqrt{\lambda_{n}-\frac{1}{4}}, & \text{if }\lambda_{n}\geq\frac{1}{4},
\\[6pt]
i\sqrt{-\lambda_{n}+\frac{1}{4}}, & \text{if }\lambda_{n}<1/4,%
\end{array}
\right.
\end{equation*}
so that $\lambda_{n}=\frac{1}{4}+r_{n}^{2}.$

Also let $\overline{m}=\max\{ k:\lambda_{k}<1/4\}.$

Let $\mathcal{G}( M) $ be the set of all closed geodesics on $M$,
and let $\mathcal{P}( M) $ be the subset of all prime closed
geodesics (that is, the closed geodesics that cannot be represented as a
non-trivial multiple of another closed geodesic). It is known that $\mathcal{G}( M)$
is a
countable set, which we can order by the lengths of its elements. Closed
geodesics correspond to hyperbolic elements of the group $\Gamma$
(that is,
the elements of $\Gamma$ with the trace outside of $[-2,2]$) up to
conjugacy of these elements. If $P\in\Gamma$ corresponds to a
geodesic $%
\gamma,$ then $\gamma$ is prime if and only if there is no $P_{0}\in
\Gamma$ such that $P=P_{0}^{k}$ for an integer $k>1.$

If $l( \gamma) $ denotes the length of the geodesic $\gamma,$ corresponding to $P \in \Gamma$,
then we set
\begin{equation*}
\left\vert\gamma\right\vert:=e^{l\left( \gamma\right) },
\end{equation*}
and note that 
\begin{equation*}
\left\vert\gamma\right\vert^{1/2}+\left\vert\gamma\right\vert
^{-1/2}=\left\vert\mathrm{Tr}P\right\vert.
\end{equation*}
We will also write $\vert\gamma\vert=N[ P] $
(meaning norm of $P$).

The Selberg trace formula relates sums over eigenvalues $\lambda_{k}$ to
sums over hyperbolic elements (geodesics) $[ P] $. Let
$h(u) $ be a function which (i) is analytic in the strip $\vert
\mathrm{Im}u\vert\leq1/2+\delta,$ (ii) is even: $h(u)=h(
-u),$ and (iii) declines sufficiently fast in the strip: $\vert
h( u) \vert=O( ( 1+\vert\mathrm{Re}%
u\vert) ^{-2-\delta}) $. 

Let $\widehat{h}( t) =\frac{1}{2\pi}\int h( u)
e^{-itu}du.$ Then the Selberg trace formula holds (cf. Theorem I.7.5 in
Hejhal \cite{hejhal76}),
%
%e19 ###
\begin{eqnarray}
\sum_{n=0}^{\infty}h\left( r_{n}\right) &=&\frac{\mu\left( F\right)
}{2\pi
}\int_{\mathbb{R}}rh\left( r\right) \tanh\left( \pi r\right) dr
\label{formula_Selberg_trace} \\
&&+\sum_{\left[ T\right] }\frac{\ln N\left[ T_{0}\right] }{N\left[
T\right]
^{1/2}-N\left[ T\right] ^{-1/2}}\widehat{h}\left( \ln N\left[ T\right]
\right), \notag
\end{eqnarray}
where the sum is over all distinct conjugacy classes of hyperbolic
elements $%
[ T] $, $[ T_{0}] $ is the primitive element for
$T,$ $%
T=T_{0}^{k},$ and $\mu( F) $ is the area of the
fundamental region of the group $\Gamma.$

It is instructive to compare this formula with formula (\ref%
{formula_Rudnick_Sarnak}). Since (\ref{formula_Selberg_trace}) resembles the explicit formulas from
number theory, it is natural to define \emph{Selberg's zeta
function} as
follows (cf. Definition II.4.1 in \cite{hejhal76}):
%
%e20 ###
\begin{equation}
Z\left( s\right) =\prod_{\gamma\in\mathcal{P}\left( M\right)
}\prod_{k=0}^{\infty}\left( 1-\left\vert\gamma\right\vert
^{-s-k}\right),%
\text{ }\mathrm{Re}s>1. \label{zeta_Selberg}
\end{equation}

It turns out that Selberg's zeta function is closely related to the
eigenvalues of the Laplace operator on $M$ (cf. Theorem II.4.10 and II.4.11
in \cite{hejhal76}).

\begin{theorem}[Hejhal-Selberg]
(a) $Z( s) $ is an entire function;
\begin{itemize}
\item[(b)] Let $\beta$ be a real number $\geq2.$ For all $s,$ the following
identity holds:
\begin{eqnarray*}
\frac{1}{2s-1}\frac{Z^{\prime}\left( s\right) }{Z\left( s\right) }
&=&\frac{%
1}{2\beta}\frac{Z^{\prime}\left( \frac{1}{2}+\beta\right) }{Z\left(
\frac{%
1}{2}+\beta\right) }+\sum_{n=0}^{\infty}\left[ \frac
{1}{r_{n}^{2}+\left( s-%
\frac{1}{2}\right) ^{2}}-\frac{1}{r_{n}^{2}+\beta^{2}}\right] \\
&&{}+\frac{\mu\left( F\right) }{2\pi}\sum_{k=0}^{\infty}\left[ \frac
{1}{%
\beta+\frac{1}{2}+k}-\frac{1}{s+k}\right].
\end{eqnarray*}
\item[(c)] $Z( s) $ has ``trivial'' zeros $s=-k,$ $k\geq1,$ with
multiplicity $( 2g-2) (
2k+1);$
\item[(d)] $s=0$ is a zero of multiplicity $2g-1;$
\item[(e)] $s=1$ is a zero of multiplicity $1;$
\item[(f)] the nontrivial zeros of $Z( s) $ are located at $\frac{1}{2}%
\pm ir_{n}.$
\end{itemize}
\end{theorem}

Since all but a finite number of eigenvalues are greater than $1/4$ hence
all but a finite number of $r_{n}$ is real and therefore the claim (f)
implies that all but a finite number of zeros of $Z( s) $ are
located on the line $\mathrm{Im}z=1/2.$

The formula in claim (b) of this theorem follows from Selberg's trace
formula and it can be thought as a functional equation for the logarithmic
derivative of $Z( s) $. In particular, it implies the functional
equation for the zeta functions itself (cf. Theorem 4.12 in \cite{hejhal76}).

\begin{theorem}[Hejhal-Selberg]
Selberg's zeta function satisfies the following functional equation:
\end{theorem}
\begin{equation*}
Z\left( s\right) =Z\left( 1-s\right) \exp\left[ \mu\left( F\right)
\int_{0}^{s-\frac{1}{2}}v\tan\left( \pi v\right) dv\right].
\end{equation*}

%s3.2 ###
\subsection{Statistics of zeros}

The number of zeta zeros in a long interval has the following asymptotic expression:
\begin{equation*}
N\left[ k:0\leq r_{k}\leq T\right] =\frac{\mu\left( F\right) }{4\pi}%
T^{2}+S\left( T\right) +E\left( T\right),
\end{equation*}
where
\begin{equation*}
S\left( T\right) =\frac{1}{\pi}\arg Z\left( \frac{1}{2}+iT\right),
\end{equation*}
and
\begin{equation*}
E\left( T\right) =O\left( 1\right) =2c\int_{0}^{T}t\left[ \tanh\left(
\pi
t\right) -1\right] dt-\left( \overline{m}+1\right).
\end{equation*}

In other words, the number of zeros in a unit interval is $\sim cT.$ In
comparison, for Riemann's zeta function we have $\sim c\log T$ zeros in the unit
interval.

It is known (cf. Theorems 8.1 and 17.1 in \cite{hejhal76}) that
\begin{equation*}
S\left( T\right) =O\left[ \frac{T}{\log T}\right],\text{ and }S\left(
T\right) =\Omega_{\pm}\left[ \left( \frac{\log T}{\log\log T}\right)
^{1/2}\right].
\end{equation*}
(Recall that the notation $f( x) =\Omega_{+}[
g(x)] $
means that $\lim\sup( \frac{f( x) }{g( x) }
) >0,$ and $f( x) =\Omega_{-}[ g(x)] $ means
that $\lim\inf( \frac{f( x) }{g( x) })
<0. $)

It was found that it is difficult to generalize the results concerning the
moments of Riemann's $S( x) $ function to the case of Selberg's zeta. Since
these results are essential for the study of statistical properties of zeta
zeros, there is a stumbling block here.

Selberg managed to resolve this problem partially for a
particular choice of the group $\Gamma.$

Let $p\geq3$ be a prime and $A$ be a quadratic non-residue modulo $p.$
Define
\begin{equation*}
\Gamma=\Gamma\left( A,p\right) =\left\{\!\left(\!\!
\begin{array}{c@{\ \ }c}
y_{0}+y_{1}\sqrt{A} & y_{2}\sqrt{p}+y_{3}\sqrt{Ap} \\
y_{2}\sqrt{p}-y_{3}\sqrt{Ap} & y_{0}-y_{1}\sqrt{A}%
\end{array}
\!\!\right)\!;\text{ }y_{0},y_{1},y_{2},y_{3}\text{ are integer.}\right\}
\end{equation*}
and call it a \emph{quaternion group.}

Let $S( t) =S^{+}( t) -S^{-}( t),$
where $%
S^{+}( t) =\max\{ 0,S( t) \} $ and $%
S^{-}( t) =\max\{ 0,\break-S( t) \}.$ Then the
following theorem holds (cf. Theorem 18.8 in \cite{hejhal76}).

\begin{theorem}[Hejhal-Selberg]
Let $\Gamma=\Gamma( A,p) $ with $p\equiv1$ $( \mathrm
{mod}%
4).$ Then (for large $T$):
\begin{equation*}
\frac{1}{T}\int_{T}^{qT}S^{+}\left( t\right) ^{2}dt\geq c_{1}\frac
{T}{\left(
\log T\right) ^{2}}
\end{equation*}
where $c_{1}$ is a positive constant that depends only on $\Gamma.$ A
similar inequality holds for $S^{-}( t).$
\end{theorem}

In order to appreciate this result note that it suggests that the average
deviation of $S( T) $ from its mean is of the
order larger than $\sqrt{T}/( \log T) $ which should be compared
with the number of zeros in the interval $[ 0,T],$ that is,
$%
cT^{2}.$ In other words, the deviation is larger than $( \mathcal
{N}%
( T) ) ^{1/4-\varepsilon}.$ To put it in prospective note
that the average deviation of the zeros of $S( T) $ for Riemann's
zeta function is of the order $( \log\log T) ^{1/2}$ which is
smaller than $\log\log\mathcal{N}( T),$ where $\mathcal{N}
( T) \sim cT\log T$ is the number of Riemann's zeros in
$[0,T].$ These situations appear to be quite different.

Moreover, recently there was some numeric and theoretical work on the
eigenvalues of the Laplace operator on manifolds $\Gamma\backslash
\mathbb{H%
}$ for arithmetic groups $\Gamma.$ First, numeric and heuristic analysis
showed that the spacings between eigenvalues resemble spacings between
points from a Poisson point process rather than spacings between eigenvalues
of a random matrix ensemble (see Bogomolny et~al.~\cite{bggs92} and
references wherein).\vadjust{\goodbreak} Next some rigorous explanations of this finding have
been given that relate it to large multiplicities of closed
geodesics with the same length. See Luo and Sarnak (\cite
{luo_sarnak94a} and
\cite{luo_sarnak94b}) and Bogomolny et al. \cite{bogomolny_leyvraz_schmit96}.

There is also some work on correlations of closed geodesics -- see Pollicott
and Sharp \cite{pollicott_sharp06}.

%s3.3 ###
\subsection{Comparison with the circle problem}

The Selberg zeta function is closely related to counting geodesics on a
space $\Gamma\backslash\mathbb{H}$, in the same way as the Riemann zeta
function is related to counting primes. It is natural to look at Laplace
eigenvalues and geodesic counting problem in a simpler situation, such
as a
compact Riemann surface of genus $1$. Such a surface can be represented as
a quotient space $\Lambda\backslash\mathbb{C}$, where $\Lambda$ is a
lattice. Consider, for concreteness, $\Lambda=[ 1,i].$ Then the
eigenvalues of the Laplace operator are $4\pi^{2}( m^{2}+n^{2}), $
where $m$ and $n$ are integer, and the number of the eigenvalues below $
t $ equals the number of integer points in the circle $t/\pi.$ Let
\begin{equation*}
r\left( n\right) =N\left\{ \left( a,b\right) \in\mathbb{Z}\times
\mathbb{Z}%
:a^{2}+b^{2}=n\right\},
\end{equation*}
and
\begin{equation*}
A\left( x\right) =\sum_{0\leq n\leq x}r\left( n\right) =\pi x+R\left(
x\right).
\end{equation*}

The function $A( x) $ can be thought as the counting function
both for eigenvalues of the Laplace operator and for closed geodesics of
bounded length.

Then by using the Poisson summation formula it is possible to derive the
following result (cf. Theorem 4.1 in \cite{hejhal76b} and Theorem 559 in
\cite{landau47}).
\begin{equation*}
\sum r\left( n\right) f\left( n\right) =\pi\sum r\left( n\right)
\int_{0}^{\infty}f\left( x\right) J_{0}\left( 2\pi\sqrt{nx}\right) dx.
\end{equation*}
Informally, if one uses this identity with the indicator function for $%
f( x) $ (which is, in fact, not allowed under the conditions of
the theorem), then one obtains the following formula (cf. formula
(4.10) in
\cite{hejhal76b} )
\begin{equation*}
R\left( x\right) =\sqrt{x}\sum_{n=1}^{\infty}\frac{r\left( n\right)
}{\sqrt{%
n}}J_{1}\left( 2\pi\sqrt{nx}\right).
\end{equation*}
Rigorous variants of this formula lead to various estimates on $R(
x),$ in particular it is known (cf. Theorems 509, 542, and 548 in
\cite{landau47}) that
\begin{equation*}
R=O\left( x^{1/3}\right) \text{ and }R=\Omega_{\pm}\left(
x^{1/4}\right),
\end{equation*}
and that
\begin{equation*}
\frac{1}{x}\int_{0}^{x}R\left( t\right) ^{2}dt=cx^{1/2}+O\left[ \left(
\log
x\right) ^{3}\right].
\end{equation*}

This suggest that the ``standard
deviation'' of $R( t) $ is $x^{1/4}.$ Similar to
the case with Laplacian eigenvalues on $\Gamma\backslash\mathbb{H}$, the
statistical behavior of eigenvalues does not resemble the behavior of random
matrix eigenvalues or Riemann's zeta zeros.\looseness=-1\vadjust{\goodbreak}

Some more details about this problem can be found in \cite{kendall48}, which
considers the question about the number of points inside a random circle.
More recent research can be found in \cite{heath-brown92}, where it is shown
that the distribution of the error term $R( x) $ converges
to a
non-Gaussian distribution as $x\rightarrow\infty$, and in \cite{bcdl93},
where this result is extended to circles with the center at a point $%
( \alpha,0)$, and it is shown that the nature of the resulting
distribution depends strongly on $\alpha.$

%s4 ###
\section{Zeta functions of dynamical systems}

Dynamical zeta functions are closely related to Selberg's zeta function which
can be
thought as a dynamical zeta function for the geodesic flow on a Riemann surface. At the same
time, there is a connection to number-theoretical zeta functions,
namely, to
the zeta functions of curves over finite fields. The main sources for this
section are reviews by Ruelle (\cite{ruelle92} and \cite{ruelle02}) and
Pollicott (\cite{pollicott03} and \cite{pollicott10b}).

%s4.1 ###
\subsection{Zetas for maps}

Let $f$ be a map of a set $M$ to itself, let the periodic orbits of $f$ be
denoted by $P,$ and let $\vert P\vert$ denote the period of the
orbit $P.$ Then, we can define the zeta of $f$ by the following formula:
%
%e21 ###
\begin{eqnarray}
\zeta\left( z\right) &=&\prod_{P}\left( 1-z^{\left\vert P\right\vert
}\right) ^{-1} \label{dynamic_zeta_definition} \\
&=&\exp\sum_{m=1}^{\infty}\frac{z^{m}}{m}\left\vert\mathrm{Fix}\text
{ }%
f^{m}\right\vert, \notag
\end{eqnarray}
where $\vert\mathrm{Fix}\text{ }f^{m}\vert$ denote the number
of fixed points of $f^{m}.$

%s4.1.1 ###
\subsubsection{Permutations}

Let $M$ be a finite set, and let $f$ be given by a permutation matrix $A.$
Then the number of fixed points of $f^{m}$ is given by \textrm{Tr}$A^{m}.$
Hence, we have
\begin{eqnarray*}
\zeta\left( z\right) &=&\exp\mathrm{Tr}\sum_{m=1}^{\infty}\frac{\left(
zA\right) ^{m}}{m} \\
&=&\exp\left( -\mathrm{Tr\log}\left( 1-zA\right) \right) \\
&=&1/\det\left( 1-zA\right),
\end{eqnarray*}
which is closely related to the characteristic polynomial of matrix
$A.$ In
particular, all poles of the zeta are on the unit circle.

%s4.1.2 ###
\subsubsection{Smooth mappings of compact manifolds}

Let $f$ be a differentiable mapping of a compact orientable smooth
manifold $%
X$ to itself. Assume that that $f$ is non-singular at all fixed
points.
Recall that the degree of $f$\vadjust{\goodbreak} at a fixed point $x$ equals to $+1$ if
the map
preserves orientation at the fixed point, and to $-1$ if it reverses the
orientation, that is, $\deg_{x}( f) :=\mathrm{sign}\det
(df_{x}-I) .$ We define the \emph{Lefschetz zeta function} as
\begin{equation*}
\zeta_{L}\left( z\right) =\exp\sum_{m=1}^{\infty}\frac{z^{m}}{m}%
\sum_{x\in\mathrm{Fix}\left( f^{m}\right) }d_{x}\left( f^{m}\right).
\end{equation*}
In this case one can use the Lefschetz fixed point formula that says:
\begin{equation*}
\sum_{x\in\mathrm{Fix}\left( f^{m}\right) }d_{x}\left( f^{m}\right)
=\sum_{i=0}^{\dim M}\left( -1\right) ^{i}\mathrm{Tr}\left( \left(
f^{m}\right) _{\ast i}:H_{i}\rightarrow H_{i}\right),
\end{equation*}
where $H_{i}$ is the $i$-th homology group of the compact manifold $M$ with
real coefficients, and $( f^{m}) _{\ast i}$ is the map
induced by
$f^{m}$ on $H_{i}.$

In particular, if $\lambda_{ij}$ are eigenvalues of $f_{\ast i},$ then we
get
\begin{eqnarray*}
\zeta_{L}\left( z\right) &=&\exp\sum_{m=1}^{\infty}\frac{1}{m}%
\sum_{i=0}^{\dim M}\left( -1\right) ^{i}\sum_{j=1}^{\dim H_{i}}\left(
z\lambda_{ij}\right)^m \\
&=&\prod_{i=0}^{\dim M}\left( \prod_{j=1}^{\dim H_{i}}\left( 1-z\lambda
_{ij}\right) ^{-1}\right) ^{\left( -1\right) ^{i}}, \\
&=&\prod_{i=0}^{\dim M} \det\left(1-z f_{\ast i} \right)^{(-1)^{i+1}}
\end{eqnarray*}
which is a rational function. If the map $f$ is a complex-analytic map of
two complex compact manifold then this calculation can be refined by using
the holomorphic Lefschetz formula that relates a sum over the fixed points
of such a map to a sum over its Dolbeault cohomology groups. This often leads
to additional information about $\lambda_{ij}.$

As an example, let $M$ be a torus $\mathbb{R}^{2}/\mathbb{Z}^{2}$ and
let $f$
be induced by a linear transformation $A\in SL_{2}( \mathbb
{Z}).$
Assume that the eigenvalues of $A$ are positive and not on the unit
circle: $%
\lambda_{1}>1>\lambda_{2}>0.$ Then $\sum_{x}\deg_{x}(
f^{m})
=\det( A^{m}-I) $. Hence, we have
\begin{eqnarray*}
\zeta_L \left( z\right) &=&\exp\sum_{m=1}^{\infty}\frac{z^{m}}{m}\det
\left(
A^{m}-I\right) \\
&=&\frac{\left( 1-z\lambda_{1}\right) \left( 1-z\lambda_{2}\right) }{%
\left( 1-z\right) \left( 1-z\lambda_{1}\lambda_{2}\right) } \\
&=&\frac{\left( 1-z\lambda_{1}\right) \left( 1-z\lambda_{2}\right) }{%
\left( 1-z\right) ^{2}}.
\end{eqnarray*}

The original dynamical zeta of continuous maps (in which fixed points are
counted without taking into account the degree of $f^{m}$) is often called
the \emph{Artin-Mazur zeta function} (see Artin-Mazur \cite{artin_mazur65}).
If this map is a diffeomorphism (a~bijection smooth in both directions) and
if it satisfies some additional conditions (hyperbolicity or Axiom A), then
it is known that this function is rational. (This~was conjectured by Smale
\cite{smale67}, and proved by Guckenheimer \cite{guckenheimer70} and
Manning~\cite{manning71}.)

%s4.1.3 ###
\subsubsection{Subshifts}

Suppose next that $A$ is an $N$-by-$N$ matrix of zeros and ones, and that
the set $M$ consists of doubly infinite sequences $\{ x_{i}
\} $
of symbols $1,\ldots,N,$ which satisfy the following criterion. A
sequence $%
\{ x_{i}\} $ belongs to $M$ if and only if $A_{x_{i}x_{i+1}}=1$
for every $i.$ In other words, the symbol $x_{i}$ determines which of the
other symbols are possible candidates for $x_{i+1}.$ The map $f$ is
simply a
shift on this set $M$: $\{ x_{i}\} \rightarrow\{
x_{i+1}\}.$ In this case, the number of fixed points of $f^{m}$
is $%
\mathrm{Tr}( A^{m}),$ and we have
%
%e22 ###
\begin{equation}
\xi\left( z\right) =1/\det\left( 1-zA\right).
\label{formula_det_subshift}
\end{equation}
%

%s4.1.4 ###
\subsubsection{Ihara's zeta function}

A basic reference for this section is a book by Terras \cite{terras11}.

Let $G$ be a finite graph. \emph{Ihara's zeta function} of $G$ is a
dynamical zeta function for a subshift associated with this graph. Namely,
orient edges of $G$ arbitrarily. Let the $2\vert E\vert$
oriented edges be denoted $e_{1},e_{2},\ldots
,e_{n},e_{n+1}=e_{1}^{-1},\ldots,e_{2n}=e_{n}^{-1}.$ The subshift
matrix is a $2n$-by-$2n$ edge adjacency matrix $W_{G}$ which is defined as follows.

\begin{definition}
The \emph{edge adjacency matrix} $W_{G}$ is a $2n$-by-$2n$ matrix with the
rows and columns corresponding to oriented edges such that its $(
i,j) $ entry equals $1$ if the terminal vertex of edge $i$ equals the
initial vertex of edge $j$ and edge $j$ is not the inverse of edge $i.$
\end{definition}

In particular, from (\ref{formula_det_subshift}) we have a determinantal
formula:
\begin{equation*}
\zeta_{G}\left( u\right) ^{-1}=\det\left( I-uW_{G}\right).
\end{equation*}

It is possible to define Ihara's zeta function more directly. The finite points of $f^{m}$ in this
example are
closed non-backtracking tailless paths of length $m,$ where a \emph
{path} is
a sequence of oriented edges such that the end of one edge equals the
beginning of the next edge. A path $( e_{1},e_{2},\ldots
,e_{m}) $
is \emph{closed} if the end of $e_{m}$ corresponds to the beginning of $
e_{1} $. It is \emph{non-backtracking} if $e_{i+1}\neq e_{i}^{-1}$ for
any $%
i,$ and it is \emph{tailless} if $e_{m}\neq e_{1}^{-1}.$

Hence, according to (\ref%
{dynamic_zeta_definition}),
%
%e23 ###
\begin{equation}
\zeta_{G}\left( u\right) =\prod_{\left[ P\right] }\left( 1-u^{l\left(
P\right) }\right) ^{-1}, \label{Ihara_zeta_definition}
\end{equation}
where the product is over all primes, that is, all equivalence classes of
primitive closed non-backtracking tailless paths. (A closed path $P$ is
\emph{primitive} if $P\neq D^{m}$ for $m\geq2$ and any other path $D;$ and
two paths are called equivalent if they can be obtained from each other
by a
cyclic permutation of edges.)\vadjust{\goodbreak}

For arbitrary regular finite graphs, this definition of Ihara's zeta
function was introduced by Sunada (\cite{sunada86} and \cite{sunada88}))
following a suggestion in the book ``Trees'' by Serre.

Ihara's zeta function has another representation as a determinant which
was first discovered by Ihara for regular graphs \cite{ihara66} and then
by Bass \cite{bass92} and Hashimoto \cite{hashimoto89} for arbitrary finite
graphs:
\begin{equation*}
\zeta_{G}\left( u\right) ^{-1}=\left( 1-u^{2}\right) ^{\left\vert
E\right\vert-\left\vert V\right\vert}\det\left(
I-A_{G}u+Q_{G}u^{2}\right),
\end{equation*}
where $A_{G}$ is the adjacency matrix of $G,$ and $Q_{G}$ is the diagonal
matrix whose $j$-th diagonal entry is ($-1$ + degree of $j$-th vertex).

If the graph is $q+1$ regular, that is, if every vertex has degree $q+1,$
then $Q_{G}$ is scalar and we can see that poles of $\zeta_{G}(
u) $ are related to the zeros of the characteristic polynomial of $A_{G},$
that is to the eigenvalues of the matrix $A_{G.}$ Precisely, the poles $u_{i}$
are related to the eigenvalues $\lambda_{i}$ by the formula:
\begin{equation*}
u_{i}=\frac{\lambda_{i}\pm\sqrt{\lambda_{i}^{2}-4d}}{2d}.
\end{equation*}
Since the eigenvalues are always real, we find that a pole $u_{i}$ is
on the
circle $\vert u\vert=1/\sqrt{d}$ if and only if the
corresponding eigenvalue is sufficiently small, $\vert\lambda
_{i}\vert\leq2\sqrt{d}.$

There are two trivial poles at $1$ and $1/q$ corresponding to the largest
eigenvalue $\lambda=d+1.$ The Riemann hypothesis for regular graphs says
that all non-trivial poles are on this circle. This is not always true, and
it holds if and only if $\vert\lambda_{1}\vert\leq2\sqrt{q},$
$\ $where $\lambda_{1}$ is the second largest in magnitude eigenvalue
of $%
A_{G}$. The graphs that satisfy this condition are often called Ramanujan graphs
following a paper by Lubotsky, Phillips, and Sarnak \cite%
{lubotzky_phillips_sarnak88}, which constructed an infinite family of such
graphs by using the Ramanujan conjecture from the theory of modular forms.

A random regular graph is approximately Ramanujan with high probability.
This means that for arbitrary $\varepsilon>0,$ the probability that $%
\vert\lambda_{1}\vert\geq2\sqrt{q}+\varepsilon$ becomes
arbitrarily small as the size of the graph grows. (This is known as Alon's
conjecture, and was proved in a lengthy paper by Friedman \cite
{friedman08}%
). The distribution properties of the largest eigenvalue are still unknown.
It is also unknown what proportion of the eigenvalues exceed the
threshold $2%
\sqrt{q}.$

In his Ph.D. thesis \cite{newland05}, Derek Newland studied spacings of
eigenvalues of random regular graphs and spacings of the Ihara zeta zeros
and found numerically that they resemble spacings in the Gaussian Orthogonal
Ensemble. See also an earlier paper by Jacobson, Miller, Rivin and Rudnick
\cite{jmrr99}.

For Ihara'z zeta function, there is an analog of Selberg's trace formula
which we formulate for the case of\ a $q+1$ regular graph.

First, note that formula (\ref{Ihara_zeta_definition}) implies
\begin{eqnarray*}
u\frac{\zeta_{G}^{\prime}}{\zeta_{G}}\left( u\right) &=&\sum_{\left[
P%
\right] }l\left( P\right) \left( u^{l\left( P\right) }+u^{2l\left(
P\right)
}+\cdots\right) \\
&=&\sum N_{m}u^{m},
\end{eqnarray*}
where $N_{m}$ be the number of all non-backtracking tailless closed
paths of
length~$m.$ Then the following explicit formula holds (cf. Proposition 25.1
in \cite{terras11}).

\begin{theorem}[Terras]
Suppose $0<a<1/q$. Assume that $h( u) $ is meromorphic in the
plane and holomorphic for $\vert u\vert>a-\varepsilon,$ $%
\varepsilon>0$ Assume that \ $h( u) =O( \vert
u\vert^{-1-\alpha}),$ $\alpha>0.$ Let $\widehat{h}(
n) :=\frac{1}{2\pi i}\int_{\vert u\vert=a}u^{n}h(
u) du$ and assume that $\widehat{h}( n) $ decays rapidly
enough. Then,
\begin{equation*}
\sum_{\rho}\rho h\left( \rho\right) =\sum_{n\geq1}N_{n}\widehat
{h}\left(
n\right),
\end{equation*}
where the sum on the left is over the poles of $\zeta_{G}( u).$
\end{theorem}

The book \cite{terras11} mentions that this formula can be used to derive
the limit law for eigenvalues of a large random regular graph (the
McKay-Kesten law) and give references to this and some other applications.

%s4.1.5 ###
\subsubsection{Frobenius maps}

Let $M$ be the set of solutions of a system of algebraic equations in $r$
variables over the algebraic closure of the finite field $\mathbb{F}_{q}$
and let $f$ be the Frobenius map $\mathrm{Frob}$: $( x_{1},\ldots
,x_{r}) \rightarrow( x_{1}^{q},\ldots,x_{r}^{q}).$ In
the case of curves, the dynamic zeta function defined by (\ref%
{dynamic_zeta_definition}) is equivalent to a number-theoretic zeta function
introduced by Artin.

Namely, let an affine curve $C$ be given by the equation $f(
X,Y)
=0$ over the finite field $\mathbb{F}_{q}.$ Let $\mathfrak{p}$ denote a
prime ideal of the field $\mathbb{F}_{q}[ X,Y] /f(
X,Y) $ and let the order of $\mathbb{F}_{q}[ X,Y] /%
\mathfrak{p}$ be denoted by $N\mathfrak{p.}$ Then, by analogy with Riemann's
zeta function we can define
\begin{equation*}
\zeta_{C}\left( s\right) =\prod_{\mathfrak{p}}\frac{1}{1-N\mathfrak
{p}^{-s}}%
.
\end{equation*}

It turns out that if
one uses the change of variable $u=q^{-s},$ then this function is equivalent to the dynamical zeta function
for the Frobenius map acting on the algebraic closure of the curve $C.$ 

Here is an example. Let $C$ be an elliptic curve, then
\begin{equation*}
\zeta_{C}\left( u\right) =\frac{1-u\left( q+1-N_{1}\right)
+u^{2}q}{\left(
1-u\right) \left( 1-uq\right) },
\end{equation*}
where $N_{1}$ is the number of the points of $C$ in $\mathbb{F}_{q}.$ One proof of
this fact is based on the Riemann-Roch theorem that allows to count the ideals
with a given norm more or less explicitly. (See \cite{milne06}.) Another
proof uses an analogy with the smooth maps of compact manifolds. It proceeds
by constructing a theory of cohomologies for algebraic varieties over finite
fields, which has a suitable Lefschetz fixed point formula. See the
book by
Silverman \cite{silverman86} for more details about this proof.

The Riemann hypothesis in this example is equivalent to the statement that
\begin{equation*}
\left\vert N_{1}-q-1\right\vert\leq2\sqrt{q},\vadjust{\goodbreak}
\end{equation*}
since this implies that the roots of the polynomial in the numerator
are on
the circle with radius $q^{-1/2}.$ For elliptic curves this was proved by
Hasse.

In \cite{weil49} Weil conjectured that the zeta function of every algebraic
variety over finite fields is rational, that it has a functional equation,
and that it satisfies the Riemann hypothesis. The proof of this general
conjecture led to an introduction of many new ideas in algebraic
geometry by
Dwork, Grothendick, and Deligne.

The statistical properties of the zeta's zeros were investigated in the case
of curves over finite fields by Katz and Sarnak \cite{katz_sarnak99}.  They
have studied zeros' distribution when the genus of the corresponding curve
grows to infinity and found that for ``most'' of the curves the local
statistics of the zeros
approach those of the eigenvalues of random matrix ensembles.

%s4.1.6 ###
\subsubsection{Maps of an interval}

For yet another example consider the map $x\rightarrow1-\mu x^{2}$ of the
interval $[ -1,1] $ to itself. For a special value of $\mu
\approx1.401155\ldots$ (the Feingenbaum value), this map has one periodic
orbit of period $2^{n}$ for every integer $n\geq0.$ Therefore, for the
dynamic zeta function we have:
\begin{equation*}
\zeta\left( z\right) =\prod_{n=0}^{\infty}\left( 1-z^{2^{n}}\right)
^{-1}=\prod_{n=0}^{\infty}\left( 1+z^{2^{n}}\right) ^{n+1}.
\end{equation*}
This $\zeta$ satisfies the functional equation $\zeta(
z^{2})
=( 1-z) \zeta( z).$ More generally, the piecewise
monotone maps of the interval $[ -1,1] $ to itself
correspond to
the \emph{Milnor-Thurston zeta functions}. See \cite{milnor_thurston88}.
Apparently, so far there have been no systematic study of the statistical
properties of their poles and zeros.

%s4.2 ###
\subsection{Zetas for flows}

If $f$ is a flow on $M,$ that is, a map $M\times\mathbb
{R}^{+}\rightarrow
M, $ then we can define the zeta function of this flow as
\begin{equation*}
\zeta\left( s\right) =\prod_{\omega}\left( 1-e^{-sl\left( \omega
\right)
}\right) ^{-1},
\end{equation*}
where $\omega$ denotes a periodic orbit of $f,$ and $l( \omega
)$ is its length. It is a more general case than the case of maps since there
is a construction (``suspension'') that
allows us to realize maps as flows (see Example 5 on p. 70 in \cite%
{pollicott_yuri98} or Section 4.3 in \cite{pollicott10b}) but not vice
versa. Unsurprisingly, it turns out that zeta functions for flows are more
difficult to investigate than zeta functions for maps.

If we imagine that prime numbers correspond to periodic orbits of a
flow and
that the length of the orbit indexed by $p$ is given by $\log p,$ then the
zeta function of the flow will coincide with Riemann's zeta function.

One particularly important example of a flow is the geodesic flow on a
smooth manifold $M.$ In the case when $M$ has a constant negative curvature,
the corresponding dynamical zeta function is closely related to the Selberg
zeta function. Namely,
\begin{equation*}
\zeta\left( s\right) =\frac{Z\left( s+1\right) }{Z\left( s\right)
}\text{
and }Z\left( s\right) =\prod_{n=0}^{\infty}\zeta\left( s+n\right) ^{-1},
\end{equation*}
where $Z( s) $ is as defined in (\ref{zeta_Selberg}) (see Remark
2.5 in \cite{pollicott10b}). Another important example is the geodesic flow
for billiards on polygons.

I am not aware about the functional equation for zeta functions that comes
from general geodesic flows (other than flows on manifolds of a constant
negative curvature). \ Also, it appears that not much is known about the
statistical properties of the distribution of zeros and poles of these zeta
functions.

On the other hand one can study the location of the pole with the maximal
real part and the results of this study give valuable information about the
distribution of closed geodesics. In this way, one can study the
distribution of closed geodesics on spaces of variable curvature (see
Corollary 6.11 in \cite{pollicott10b}).

%s5 ###
\section{Conclusion}

We considered the statistical properties of various zeta functions. For the
zeros of number-theoretical zeta functions, the main observation is that
they satisfies many properties which are true for eigenvalues of random
matrices. The main outstanding problem (besides the Riemann hypothesis) is
to push this similarity to its natural limits and, in particular, show that
the Montgomery conjecture about the correlations of the Riemann zeros is true.

Next, we observed that for some of groups $\Gamma\subset SL_{2}(
\mathbb{Z%
}),$ the statistical properties of Selberg's zeta zeros are
different from those of the random matrix eigenvalues. The exact description
of these properties is not known.

The statistical properties of the zeros of dynamical zeta functions are not
investigated in many cases. Two notable exceptions are the zeta functions of curves
over finite fields and Ihara's zeta functions. However, even in this case
there are many unsolved problems. For example, it is not known whether the
distribution of local statistics for the zeros of Ihara's zeta functions
coincide with the corresponding distribution for the random matrix
eigenvalues.

\bibliographystyle{plain}

\end{document}